\documentclass[11pt]{amsart}
\usepackage{amsxtra}
\usepackage{amssymb}
\usepackage{mathtools}
\theoremstyle{plain}
\newcommand{\cleqn}{\setcounter{equation}{0}}
\newcommand{\clth}{\setcounter{theorem}{0}}
\newcommand {\sectionnew}[1]{\section{#1}\cleqn\clth}
\newcommand{\nn}{\hfill\nonumber}
\newtheorem{theorem}{Theorem}[section]
\newtheorem{lemma}[theorem]{Lemma}
\newtheorem{definition-theorem}[theorem]{Definition-Theorem}
\newtheorem{proposition}[theorem]{Proposition}
\newtheorem{corollary}[theorem]{Corollary}
\newtheorem{definition}[theorem]{Definition}
\newtheorem{example}[theorem]{Example}
\newtheorem{remark}[theorem]{Remark}
\newtheorem{conjecture}[theorem]{Conjecture}

\newcommand \bth[1] { \begin{theorem}\label{t#1} }
\newcommand \ble[1] { \begin{lemma}\label{l#1} }

\newcommand \bpr[1] { \begin{proposition}\label{p#1} }
\newcommand \bco[1] { \begin{corollary}\label{c#1} }
\newcommand \bde[1] { \begin{definition}\label{d#1}\rm }
\newcommand \bex[1] { \begin{example}\label{e#1}\rm }
\newcommand \bre[1] { \begin{remark}\label{r#1}\rm }
\newcommand \bcj[1] { \begin{conjecture}\label{j#1}\rm }

\renewcommand {\eth} { \end{theorem} }
\newcommand {\ele} { \end{lemma} }

\newcommand {\epr} { \end{proposition} }
\newcommand {\eco} { \end{corollary} }
\newcommand {\ede} { \end{definition} }
\newcommand {\eex} { \end{example} }
\newcommand {\ere} { \end{remark} }
\newcommand {\ecj} { \end{conjecture} }

\newcommand {\enota} { \end{notation} }
\newcommand \thref[1]{Theorem \ref{t#1}}
\newcommand \leref[1]{Lemma \ref{l#1}}
\newcommand \prref[1]{Proposition \ref{p#1}}

\newcommand \deref[1]{Definition \ref{d#1}}

\newcommand \reref[1]{Remark \ref{r#1}}
\newcommand \lb[1]{\label{#1}}
\def \Cset {{\mathbb C}}
\def \KK {{\mathbb K}}
\def \Zset {{\mathbb Z}}
\def \Nset {{\mathbb N}}
\def \Qset {{\mathbb Q}}

\def \QQ {{\mathcal{Q}}}
\def \PP {{\mathcal{P}}}

\def \VV {{\mathcal{V}}}
\def \UU {{\mathcal{U}}}

\def \SS {{\mathcal{S}}}
\def \LL {{\mathcal{L}}}

\def \XX {{\mathcal{X}}}
\def \YY {{\mathcal{Y}}}

\def \De {\Delta}   % Greek letters
\def \de {\delta}
\def \al {\alpha}
\def \be {\beta}
\def \vpi {\varpi}
\def \la {\lambda}

\def \om {\omega}

\def \ga {\gamma}
\def \de {\delta}
\def \Ga {\Gamma}

\def \sig {\sigma}

\def \ep {\epsilon}
\def \sig{\sigma}
\def \mt  {\mapsto}
\def \lra {\longrightarrow}

\def \hra {\hookrightarrow}

\def \ci  {\circ}           % duals

\def \rcor {\rangle}
\def \lcor {\langle}
\def \del {\partial}

\def \ol {\overline}
\def \wt {\widetilde}
\def \wh {\widehat}

\def \id { {\mathrm{id}} }
\def \st { {\mathrm{st}} }

\def \coeff { {\mathrm{coeff}} } 
\def \sign { {\mathrm{sign}} }
\def \Ad { {\mathrm{Ad}} }

\def \nr { {\mathrm{nrf}} }
\def \Lie { {\mathrm{Lie \,}} }

\def \g  {\mathfrak{g}}   % Lie algebra letters

\def \n  {\mathfrak{n}}

\DeclareMathOperator \Aut { {\mathrm{Aut}} }

\DeclareMathOperator \Der {{\mathrm{Der}}}

\DeclareMathOperator \tr { {\mathrm{tr}} }

\DeclareMathOperator \gr  { {\mathrm{gr}} }

\newcommand \Spec { {\mathrm{Spec}} }
\begin{document}
%%%%%%%%%%%%%%%%%%%%%%%%%%%%%%%%%%%%%%%%%%%%%%%%%%%%%%%%%%%%%%%%%%%%%%%%%%%
%%%%%%%%%%%%%%%%%%%%%%    Title    %%%%%%%%%%%%%%%%%%%%%%%%%%%%%%%%%%%%%%%%
\title[Noncommutative discriminants via Poisson primes]
{Noncommutative discriminants via Poisson primes}
\author[Bach Nguyen]{Bach Nguyen}
\address{
Department of Mathematics \\
Louisiana State University \\
Baton Rouge, LA 70803 \\
U.S.A.
}
\email{bnguy38@lsu.edu}
\author[Kurt Trampel]{Kurt Trampel}
\email{ktramp2@lsu.edu}
\author[Milen Yakimov]{Milen Yakimov}
\email{yakimov@math.lsu.edu}
\thanks{The research of K.T. has been supported by a VIGRE fellowship through the NSF grant DMS-0739382 and that of M.Y. by the NSF grant 
DMS-1303038 and Louisiana Board of Regents grant Pfund-403.}
%\date{}
\begin{abstract} We present a general method for computing discriminants of noncommutative algebras. It builds a
connection with Poisson geometry and expresses the discriminants as products of Poisson primes. 
The method is applicable to algebras obtained by specialization from families, such as quantum algebras at roots of unity. 
It is illustrated with the specializations of the algebras of quantum matrices at roots of unity
and more generally all quantum Schubert cell algebras. 
%The method
%is applicable to many classes of quantum algebras at roots of unity, such as the specializations of 
%the algebras of quantum matrices and more generally all quantum Schubert cell algebras. 
\end{abstract}
\maketitle
%%%%%%%%%%%%%%%%%%%%   Introduction   %%%%%%%%%%%%%%%%%%%%%%%%%%%%%%%%%%%%%%%%
\sectionnew{Introduction}
\lb{Intro}
%\subsection{}
%\label{1.1}
\subsection{}
\label{1.1}
In the commutative setting, the notion of discriminant plays a fundamental role in algebraic number theory, algebraic geometry and combinatorics.
We refer the reader to the book of Gelfand, Kapranov, and Zelevinsky \cite{GKZ} for extensive background. 
The notion also has an analog for noncommutative algebras that has been used in a key way in the study of
orders and lattices in central simple algebras, see Reiner's book \cite{Re}. Recently, Zhang, Bell, Ceken, Palmieri and Wang
found many additional applications of noncommutative discriminants
in the study of automorphism groups of PI algebras \cite{CPWZ1,CPWZ2} and related problems, such as the 
isomorphism problem and the Zariski cancellation problem \cite{BZ}. Finally, discriminants also have tight relations 
to the representation theory of the algebra in question.
 
However, the computation of the discriminant of a noncommutative algebra turns out to be a rather challenging problem. It has been solved for 
very few families of algebras: quasipolynomial rings and two cases of quantum Weyl algebras
 \cite{CPWZ1,CPWZ2,CYZ}. For example, even the discriminants of the algebras of quantum $2 \times 2$ matrices at roots of unity
are presently unknown. Furthermore, it was observed in \cite{CPWZ1,CPWZ2,CYZ} that in the known cases the discriminants have 
an elaborate product form, but the meaning of the individual terms remained unclear. The problem of computing noncommutative 
discriminants has so far been attacked with techniques within noncommutative algebra.

In this paper we present a general method for the computation of the discriminants of algebras obtained as specializations. 
It establishes a connection between discriminants of noncommutative algebras and Poisson geometry.
The method is applicable 
to broad classes of algebras, such as specializations of quantum algebras at roots of unity. The algebras of quantum matrices of any size 
at roots of unity are treated as special cases.

 The discriminant of an algebra $R_\ep$ obtained by specialization is an 
element of a Poisson algebra $A$ sitting inside the center of $R_\ep$. We prove that the discriminant is a Poisson normal element of $A$ 
and that, under conditions satisfied very frequently, it is a product of Poisson prime elements of $A$. This explains the intrinsic nature of the terms 
in the product formulas for discriminants.
The Poisson prime elements are closely related to Poisson hypersurfaces in the Poisson variety $\Spec A$.

Thus, the problem of computing noncommutative discriminants becomes one about the interplay between the PI algebra in question, $R_\ep$, and the 
symplectic foliation of the Poisson variety $\Spec A$. The precise form of the discriminant is determined by studying 
its evolution under the hamiltonian flows on $\Spec A$, as described in the next subsection.
\subsection{}
\label{1.2} In the remaining part of the introduction we formulate the main results of the paper. 
Let $R$ be an algebra over $\KK[q^{\pm1}]$ for a field $\KK$ 
of characteristic 0. Let $\ep \in \KK^\times$ be such that $q - \ep \in R$ is regular (i.e., not a zero divisor).
The specialization of $R$ at $\ep$ is the $\KK$-algebra $R_\ep: = R/(q -\ep)R$.
Denote the projection $\sig \colon R \to R_\ep$. The center $Z(R_\ep)$ has a canonical structure of Poisson algebra: 
\[
\{ \sig(x_1), \sig(x_2) \} := \sig \left( \frac{x_1 x_2 -x_2 x_1}{q-\ep} \right), \quad x_i \in \sig^{-1}(Z(R_\ep)).
\]
Let $C_\ep$ be a Poisson subalgebra of $Z(R_\ep)$ such that $R_\ep$ is a free $C_\ep$-module of finite rank.
A $C_\ep$-basis $\YY:=\{y_j \mid 1 \leq j \leq N \}$ gives rise to an embedding $R_\ep \hra M_N(C_\ep)$. 
(Here and below $M_n$ refers to square matrices and $M_{m,n}$ to rectangular matrices of those sizes.)
The composition of this embedding with the standard trace map leads to a trace map $\tr \colon R_\ep \to C_\ep$ which is independent 
on the choice of basis. The discriminant of $R_\ep$ over $C_\ep$ is defined
by
\[
d(R_\ep/C_\ep) := \det \big(  [\tr(y_i y_j)]_{i,j=1}^N  \big) \in C_\ep,
\]
see \S \ref{2.1} for details on trace maps and discriminants. The notation $[c_{ij}]_{i,j=1}^N$ refers to the 
square matrix with entries $c_{ij}$.
The element $d(R_\ep/C_\ep)$ is defined up to associates, because changing the basis $\YY$ multiplies the discriminant by an element of $C_\ep^\times$. 
(As usual, the units of an algebra $A$ will be denoted by $A^\times$; $a, b \in A$  are called associates, denoted by $a =_{A^\times} b$, 
if $a = u b$ for some $u \in A^\times$.)

An element $a$ of a Poisson algebra $A$ is called Poisson normal if the principal ideal $(a)$ is Poisson (i.e., is closed under the Poisson bracket). 
An element $p \in A$ is called Poisson prime if $(p)$ is a prime ideal which is also Poisson. The latter happens precisely when the zero locus of $p$ is 
a union of symplectic leaves of the Poisson variety $\Spec A$.  Background material on these notions and the notion of noetherian 
Poisson unique factorization domain is contained in \S \ref{2.2}. 
\medskip
\\
\noindent
{\bf{Theorem A.}} {\em{Let $R$ be a $\KK[q^{\pm1}]$-algebra where $\KK$ is a field of characteristic 0 and $\ep \in \KK^\times$
be such that $q - \ep \in R$ is regular.
Assume that $R_\ep := R/(q-\ep)R$ is a free module of finite rank over a Poisson
subalgebra $C_\ep$ of its center. 

(i) Then $d(R_\ep/C_\ep)$ is a Poisson normal element of $(C_\ep, \{.,.\})$.

(ii) Assume, in addition, that $C_\ep$ is a unique factorization domain as a commutative algebra or a noetherian Poisson unique factorization domain. 
Then, $d(R_\ep/C_\ep)=0$ or
\[
d(R_\ep/C_\ep) =_{C_\ep^\times} \prod_{i=1}^m p_i
\]
for some (not necessarily distinct) Poisson prime elements $p_1, \ldots, p_m \in C_\ep$.
}}
\medskip

In concrete situations $C_\ep$ is very close to being the full center of $R_\ep$ and, as a consequence, the
Poisson bracket on $C_\ep$ is very nontrivial. Because of this, the collection of Poisson primes of $C_\ep$ is only a very small subset 
of the set of all primes of the commutative algebra $C_\ep$. From this perspective, Theorem A is a rigidity theorem for fully noncommutative 
extensions (such for which $C_\ep$ is almost the full center of $R_\ep$) in the sense that their discriminants are products 
of Poisson primes of $C_\ep$ which altogether form a discrete and often a finite set.

Since Theorem A places a strong restriction on the possible form of the discriminant 
$d(R_\ep/C_\ep)$, the latter can be fully determined using a recipe of four methods described in \S \ref{3.3} (a combination of 
algebraic methods, such as filtration arguments, with geometric methods, such as evolution under hamiltonian flows). The Poisson primes of $C_\ep$ 
are determined from the symplectic foliation of the Poisson variety $\Spec C_\ep$. The symplectic foliations of Poisson varieties appearing in the theory 
of quantum groups are well understood in great generality \cite{EL1,EL2,LY,Y-Duke}.

From another perspective, Theorem A builds a bridge between discriminants of noncommutative algebras and algebraic number theory. 
In the latter situation, discriminants of number fields $\KK$ are expressed as products of 
those prime numbers $p \in \Zset$ that ramify in $\KK$. In our case we even have full control of the powers 
of the (Poisson) prime elements that appear in the factorization of the discriminant; these powers are determined by 
using hamiltonian flows on the basis of \prref{Poisson-brack}, see the discussion in \S \ref{3.3} (3). 

As a simple example, we first consider the quantum Weyl algebra $A_\ep$ with generators $x_1,x_2$ subject to the relation $x_1 x_2 = \ep x_2 x_1 + 1$ 
where $\ep \in \KK$ is a primitive $l$-th root of unity. The
induced Poisson structure on its center $Z(A_\ep) \cong \KK[z_1,z_2]$ (where $z_i := x_i^l$) is 
\[
\{ z_1, z_2 \} = l^2 \ep^{-1} (z_1 z_2-t) \quad \mbox{where} \quad
t := (1-\ep)^{-l}.
\]
It has only one Poisson prime $z_1 z_2 - t$, leading to the discriminant formula
\[
d(A_\ep, Z(A_\ep)) =_{\KK^\times} (z_1 z_2 -t)^{l(l-1)},
\]
see \S \ref{3.4}. This recovers a formula of Chan, Young and Zhang \cite{CYZ} which was obtained in a more complicated fashion.

Let $R_q[M_n]$ be the algebra of square quantum matrices with generators $x_{ij}$ and relations \eqref{qmatr}. 
The specialization $R_\ep[M_n]$ at a primitive $l$-th root of unity $\ep \in \KK$ is free over the central subalgebra $C_\ep[M_n]$, 
generated by the images of $x_{ij}^l$ under the projection map $R_q[M_n] \to R_\ep[M_n]$. Denote those 
images by $z_{ij}$.

For $i \leq j \in \Zset$, set $[i,j] :=\{i, \ldots, j \}$.
For two sets $I,J \subseteq [1,n]$, denote by $\De_{I;J}$ the minor in $z_{ij}$ 
corresponding to the rows in the set $I$ and the columns in $J$. 
\medskip
\\
\noindent
{\bf{Theorem B.}} {\em{For all fields $\KK$ of characteristic 0, odd integers $l>2$, primitive $l$-th roots of unity $\ep \in \KK$,
and positive integers $n$, the discriminant of the specialization of the algebra of square quantum matrices $R_\ep[M_n]$ is given by
\[
d(R_\ep[M_n]/C_\ep[M_n]) =_{\KK^\times} \prod_{k=1}^n \De_{[n-k+1, n]; [1, k]}^L 
\prod_{j=1}^{n-1}  \De_{[1, j] ; [n-j+1, n]}^L
\] 
where $L:= l^{n^2-1}(l-1)$.}}
\medskip

For all simple Lie algebras $\g$ and Weyl group elements $w \in W$, Lusztig \cite{L} and De Concini, Kac, and Procesi \cite{DKP2}
defined a subalgebra $\UU^-[w]$ of the quantized universal enveloping algebra $\UU_q(\g)$ which is a 
deformation of $\UU(\n_- \cap w(\n_+))$ for the nilradicals $\n_\pm$ of a pair of opposite Borel subalgebras of $\g$. 
The specialization at a primitive $l$-th root of unity $\ep \in \KK$ of the corresponding nonrestricted rational form 
will be denoted by $\UU^-_\ep[w]$. The central subalgebra of the latter generated by the images of the $l$-th powers of 
Lusztig's root vectors will be denoted by $C^-_\ep[w]$. One can canonically identify
$C^-_\ep[w]$ with the coordinate ring of the Schubert cell $B_+ w \cdot B_+ \subset G/B_+$ in the full flag variety of the 
split, connected, simply connected algebraic $\KK$-group $G$ associated to $\g$. The Schubert cell $B_+ w \cdot B_+$
is isomorphic to the unipotent subgroup $U_+ \cap w (U_-)$ of $G$ where $U_\pm$ are the unipotent radicals of $B_\pm$.
By restriction, generalized minors $\De_{\la, w \la}$ on $G$ are identified with elements of $C^-_\ep[w]$ where $\la$ is a
dominant weight.
\medskip
\\
\noindent
{\bf{Theorem C.}} {\em{Let $\g$ be a simple Lie algebra, $w$ a Weyl group element, and $l>2$ an odd integer
which is $\neq 3$ in the case of $G_2$. Let $\KK$ be a field of characteristic 0 containing 
a primitive $l$-th root of unity $\ep$. Then
\[
d(\UU^-_\ep[w]/C^-_\ep[w]) =_{\KK^\times}  \prod_{i \in \SS(w)}\De_{\vpi_i, w \vpi_i}^L
\]  
where $L: = l^{N-1}(l-1)$, $\vpi_i$ are the fundamental weights of $\g$ and $\SS(w)$ is the support of $w$.
}}

Theorem A is proved in Section \ref{Pprop-Disc}. Two companion theorems for $n$-discriminant ideals 
without freeness assumptions are obtained in \S \ref{3.5}. These results are extended to the setting 
of Poisson orders introduced by Brown and Gordon \cite{BrGo} in \S \ref{3.6}. 
Section \ref{disrc-qSchubert} contains the proof of Theorem C. It relies 
on the results of De Concini, Kac and Procesi \cite{DKP1,DKP1b,DKP2} on specializations of the 
nonrestricted rational forms of $\UU_q(\g)$ and $\UU^-[w]$ at root of unity, as well as on results on 
Poisson algebraic groups and Poisson homogeneous spaces. Theorem B for quantum matrices is 
proved in Section \ref{q-matr}. The proofs are completely elementary and do not use the results in \cite{DKP1,DKP1b,DKP2}.
This section was written so the proofs for the algebras of quantum matrices can be read without any Lie theoretic 
background. All algebras of quantum matrices $R_q[M_{m,n}]$
arise as special cases of the quantum Schubert cell algebras for $\g = {\mathfrak{sl}}_{m+n}$ and 
particular Weyl group elements. The relation between the corresponding discriminant formulas is described in \S \ref{5.5}.

The results in Theorem C concern quantum Schubert cell algebras for finite dimensional simple Lie algebras as opposed to such for 
symmetrizable Kac--Moody algebras because of the use of \cite{DKP1,DKP1b} on specializations of $\UU_q(\g)$. In the PhD theses 
of the first two named authors these results will be extended to the Kac--Moody case. This is done by 
working directly with the algebras $\UU_\ep^-[w]$ which are orders in central simple algebras while the specializations of $\UU_q(\g)$ are not 
finite over their centers for general symmetrizable Kac--Moody algebras $\g$.

In forthcoming papers we apply the techniques in the paper to derive formulas for the discriminants 
of other important classes of algebras: the specializations at roots of unity of the quantized universal enveloping algebras 
of simple Lie algebras, the corresponding quantum function algebras, and general quantum Weyl algebras.
%  and the cases of Sklyanin algebras which are PI.

The quantum Schubert cell algebras $\UU^-[w]$ have quantum cluster algebra structures \cite{GLS2,GY-big} and $C^-_\ep[w]$ are isomorphic to the 
corresponding classical cluster algebra structures on unipotent cells \cite{GLS2,GY}. In the framework of cluster algebras, Theorems B and C
prove that the discriminant $d(\UU^-_\ep[w]/C^-_\ep[w])$ is precisely the product of all frozen variables (raised to the same powers). We expect that 
the discriminants of the specializations of quantum cluster algebras \cite{BeZe} to primitive $l$-th roots of unity over the commutative subalgebras 
of $l$-th powers of cluster variables will have deep Poisson geometric properties in terms of the 
Gekhtman--Shapiro--Vainshtein Poisson structure \cite{GSV} on the corresponding cluster variety.

All algebras in the paper will be assumed to be unital. For a Poisson structure $\pi$ on an algebraic variety $V$ we will denote 
by $\{.,.\}_\pi$ the Poisson bracket on the coordinate ring of $V$.
\medskip
\\
\noindent
{\bf Acknowledgements.} We would like to thank Jianmin Chen, Ken Goodearl, Tom Lenagan and James Zhang for their very helpful comments and for 
pointing out a mistake in the first proof of \thref{last}.  We are thankful to the referee for the valuable suggestions which helped up to improve the exposition.
%%%%%%%%%%%%%%%%%%%%%%%%%%%%%%%%%%%%%%
\sectionnew{Discriminants and Poisson primes}
This section contains background material on discriminants and Poisson normal and prime elements.  
The section also contains new results that will be used later in the paper.
\lb{Discr-PUFD}
\subsection{}
\label{2.1}
An algebra with trace is an algebra $R$ with a linear map $\tr \colon R \to R$ such that 
\begin{enumerate}
\item $\tr(x y) = \tr (yx)$, 
\item $\tr(x) y = y \tr (x)$,
\item $\tr( \tr(x) y ) = \tr(x) \tr(y)$ 
\end{enumerate}
for all $x, y \in R$. Following Procesi \cite{P2}, an affine $\KK$-algebra with trace is said to be a 
$d$-Cayley--Hamilton algebra if it satisfies the $d$-th Cayley--Hamilton identity 
\[
\chi_{d,x}(x) =0 \; \; \mbox{for all}  \; \; x \in R 
\]
and $\tr(1) =d$.
If $C$ is the central subalgebra of $R$ generated by $tr(y)$ for $y \in R$, the  $d$-th characteristic polynomial $\chi_{d,x}(t) \in C[t]$ 
of $x \in R$ is defined to be
\[
\chi_{d,x}(t) := t^d - \wh{c}_1(x) t^{d-1} + \cdots + (-1)^d \wh{c}_d (x)
\]
where the functions $\wh{c}_i : R \to C$ are defined as follows. We express the $k$-th elementary
symmetric function $\sig_k$ in $d$ indeterminates $\la_1, \ldots, \la_d$ in terms of the Newton power sum functions
$\psi_j:=\la_1^j + \cdots + \la_d^j$ as
\[
\sig_k = p_k( \psi_1, \ldots, \psi_k) \quad \mbox{for} \quad p_k(t_1, \ldots, t_k) \in \Zset[(k!)^{-1}][t_1, \ldots, t_k],
\]
and then set $\wh{c}_k(x) := p_k( \tr(x), \ldots, \tr(x^k))$.

A geometric framework for the study of the representations of Cayley--Hamilton algebras was developed by Procesi \cite{P1,P2}.
We refer the reader to the detailed and very instructive expositions of De Concini--Procesi \cite[Section 4]{DP} and Le Bruyn \cite[Chapter 1]{Le}.
Orders in central simple algebras give rise to Cayley--Hamilton algebras, and quantized universal enveloping algebras at roots 
of unity belong to this class of algebras \cite{DKP1,DKP1b,DKP2}.
De Concini, Kac and Procesi \cite{DKP1,DKP1b,DKP2} studied the representation theory of quantum algebras at roots of unity 
in this framework with the additional use of Poisson geometry.

Returning to traces, slightly more generally, one can consider trace maps $\tr \colon R \to F$ where 
$F$ is a commutative algebra which is an extension of a subalgebra $C \subset Z(R)$. Such traces are required to be $C$-linear and cyclic.
They arise from embeddings of $R$ into $F$-algebras. 

\bde{discr} \cite{Re}
(i) The {\em{discriminant of the set}} $\YY := \{ y_1, \ldots, y_N \} \subset R$ is defined to be
\[
d_N(\YY : \tr) := \det \left(  [\tr(y_i y_j)]_{i,j=1}^N  \right) \in F.
\]

(ii) The {\em{$N$-discriminant ideal}} $D_N(R/C)$ is the $C$-submodule of $F$ generated by 
\[
d_N(\YY : \tr) \quad \mbox{for the $N$-element subsets} \; \; \YY \subset R. 
\]

(iii) If $R$ is a free, rank $N$ module over a subalgebra $C$ of the center $Z(R)$, then the embedding $R \hra M_N(C)$ gives rise 
to a trace map $\tr : R \to C$ which is independent on the choice of basis. The {\em{discriminant}} of $R$ over $C$ is defined by 
\[
d(R/C) :=_{C^\times} d_N(\YY : \tr)
\] 
where $\YY=\{ y_1, \ldots, y_N\}$ is a $C$-basis of $R$. 
\ede
\noindent
The discriminant in part (iii) is well defined since for any other 
$C$-basis $\XX:=\{ x_1, \ldots, x_N\}$ of $R$, 
\begin{equation}
\label{discr-bas}
d_N(\XX : \tr) = \det(b)^2 d_N(\YY : \tr)
\end{equation}
where $b:=(b_{ij}) \in M_N(C) $ is the change of bases matrix given by $x_i = \sum_j b_{ij} y_j$, see \cite[Exercise 4.13]{Re}.
Our notation in \deref{discr} (iii) differs from the one in \cite{CPWZ1}.

\bpr{invar} Assume that $R$ is an algebra with trace $\tr \colon R \to C \subset Z(R)$ which 
is a free module over $C \subset Z(R)$ of rank $N$.

(i) If $\partial$ is a derivation of $R$ such that $\tr ( \partial x) = \partial \tr (x)$ for all $x \in R$, then 
\[
\del d_N(\YY : \tr) = 2 \tr(b) d_N(\YY : \tr)
\]
for any $C$-basis $\YY := \{ y_1, \ldots, y_N \}$ of $R$,
where $b =(b_{ij}) \in M_N(C)$ is the matrix with entries given by $\del y_i = \sum_j b_{ij} y_j$.

(ii) Let $\tr \colon R \to C$ be the canonical trace map from the embedding $R \subset M_N(C)$
associated to any $C$-basis of $R$.
Then every derivation $\partial$ of $R$, satisfying $\partial(C) \subseteq C$,  
has the property 
$\tr ( \partial x) = \partial \tr (x)$, $\forall x \in R$.
\epr
\begin{proof} (i) The discriminant $d_N(\YY : \tr)$ is the same when it is computed for the pair of algebras 
$(R,C)$ and $(R[t]/(t^2), C[t]/(t^2))$. From the second interpretation and \eqref{discr-bas} we get 
\[
(1 + t \partial) d_N(\YY : \tr) = d_N((1+t\partial)\YY : \tr) = \det(I_N + t b)^2 d_N(\YY : \tr)
\]
where $(1 + t \partial) \YY:= \{ (1+t \partial) y_1, \ldots, (1+ t \partial) y_N \}$. The first part now follows from the fact 
that $d_N(\YY : \tr) \in C$ by comparing the coefficients of $t$.

(ii) Given a basis of $R$ over $C$, from the embedding $R[t]/(t^2) \hra M_N(C[t]/(t^2))$ we obtain
$(1+ t \partial) \tr (x) = \tr((1+ t \partial) x)$ for all $x \in R$ where the traces are computed 
in $R[t]/(t^2)$. This implies the statement since $\tr(x), \tr (\partial x) \in C$.
\end{proof}
The second part is valid in much greater generality for orders in central simple algebras \cite[Ch. 9-10]{Re}, but we will 
not need this here.
%%%%%%%%%%%%%%
\subsection{}
\label{2.2} Let $(A, \{.,.\})$ be a Poisson algebra over a base field of characteristic 0. It is called noetherian 
if it is noetherian considered as a commutative algebra.
\bde{prime-norm} (i) An element $a \in A$ is called {\em{Poisson normal}} if for every $x \in A$,
\[
\{a, x \} = ay \quad \mbox{for some $y \in A$}.  
\]
If $A$ is an integral domain as a (commutative) algebra, this is equivalent to saying that
\[
\{a, x \} = a \partial(x) \quad \mbox{for some Poisson derivation $\partial$ of $A$.}  
\]

(ii) Assume that $A$ is an integral domain as an algebra. An element $p \in A$ is called {\em{Poisson prime}} 
if it is a prime element of the algebra which is normal in the Poisson sense. 
\ede
\bre{prime-norm-def} (i) An element $a \in A$ is normal if and only if the principal ideal $(a)$ is Poisson.

(ii) An element $p \in A$ is Poisson prime if and only if the ideal $(p)$ is nonzero, prime and Poisson.

(iii) Assume that the base field is $\Cset$ and $\Spec A$ is smooth. View the elements of $A$ as regular functions on the Poisson variety $\Spec A$.
A prime element $p \in A$ is Poisson prime if and only if its zero locus $\VV(p)$ is a union of symplectic leaves of $\Spec A$.

This can be proved as follows. If $\VV(p)$ is a union of symplectic leaves of $\Spec A$, then for all $g \in A$, $\{p, g\}$ vanishes on the smooth locus of $p$. 
Thus, $\{p, g\}$ belongs to $(p)$.  In the opposite direction, assume that $(p)$ is Poisson. If $\LL$ is a symplectic leaf of $\Spec A$ 
such that $\LL \cap \VV(p) \neq \varnothing$ and $\LL \not\subseteq \VV(p)$, then for every smooth point $m \in \LL \cap \VV(p) \subsetneq \LL$ 
there will exist $g \in A$ such that $\{p, g\}(m) \neq 0$. This would contradict the assumption that $(p)$ is a Poisson ideal.
\ere
\bde{PUFD} A noetherian Poisson algebra $A$ is called a Poisson unique factorization domain if it is an integral domain
as an algebra and every  non-zero Poisson prime ideal of $A$ contains a Poisson prime element.
\ede
\bpr{fact} Let $A$ be a Poisson algebra over a field of characteristic 0, satisfying one of the following 2 conditions:
\begin{itemize}
\item $A$ is a unique factorization domain as a commutative algebra or
\item $A$ is a noetherian Poisson unique factorization domain. 
\end{itemize}
Then every non-zero, non-unit Poisson normal element $a \in A$ has a unique factorization 
of the form
\[
a = \prod_{i=1}^m p_i
\]
for some set of (not necessarily distinct) Poisson prime element $p_1, \ldots, p_m \in A$. The uniqueness is up to taking associates
and permutations.
\epr
The case of noetherian Poisson UFDs is analogous to the unique factorization property of normal elements 
in (noncommutative) noetherian UFDs proved by Chatters \cite[Proposition 2.1]{Cha}, see also \cite[Proposition 2.1]{GY}. 
The case when $A$ is a UFD as a commutative algebra follows from the following lemma. 
\ble{ClassUFD} Assume that $A$ is a Poisson algebra over a field of characteristic 0 which is a unique factorization domain as a commutative algebra. 
If $a \in A$ is a Poisson normal element and $p \in A$ is a prime element such that $p \mid a$, 
then $p$ is a Poisson prime element.
\ele
\begin{proof} Let $a = p^k b$ for some $b \in A$ such that $p \nmid b$. For every $x \in A$, there exists $y \in A$ such that $\{ a , x \} = a y$. Then 
\[
k \{p, x\} p^{k-1} b + p^k \{b, x \} = p^k b y.
\] 
Since the base field has characteristic 0, we have $p^k \mid \{p, x\} p^{k-1} b$ for every $x \in A$, and so $p \mid \{p, x\}$.
\end{proof}
%%%%%%%%%%%%%%%%%%%%%%%%%%%%%%%%%%%%%%
\sectionnew{General theorems on noncommutative discriminants}
\lb{Pprop-Disc}
In this section we prove two general theorems on discriminants and discriminant ideals
of algebras obtained as specializations. We also give a recipe about computing discriminants from the 
first theorem. In the last subsection we obtain extensions of these results to the setting of Poisson orders introduced by Brown and 
Gordon \cite{BrGo}. In order to keep the exposition more transparent we first prove the results in the more common 
setting of specializations, and then extend them to the setting of Poisson orders.
\subsection{}
\label{3.1}
Let $R$ be an algebra over $\KK[q^{\pm1}]$. For $\ep \in \KK^\times$, one defines the specialization of $R$ at $\ep$ which is
the $\KK$-algebra $R_\ep: = R/(q -\ep)R$. Denote the canonical projection $\sig \colon R \to R_\ep$.
Assume that $q - \ep \in R$ is a regular element. 
The center $Z(R_\ep)$ has a canonical structure of Poisson algebra defined as follows. For $z_1, z_2 \in Z(R_\ep)$, choose $x_i \in \sig^{-1}(z_i)$ and set 
\begin{equation}
\label{Poisson}
\{ z_1, z_2 \} := \sig \left( \frac{x_1 x_2 -x_2 x_1}{q-\ep} \right).
\end{equation}
This definition does not depend on the choice of $x_1$ and $x_2$. Indeed, if $\sig(x_1) = \sig(x'_1)$, then $x_1 - x'_1 = (q-\ep) x$ for some 
$x \in R$ and
\[
\sig \left( \frac{[x_1, x_2]}{q-\ep} \right) - \sig \left( \frac{[x'_1, x_2]}{q-\ep} \right) = \sig([x, x_2]) =0 
\]
because $\sig(x_2) \in Z(R_\ep)$. Similarly one checks that $\{z_1, z_2 \} \in Z(R_\ep)$ for $z_i \in Z(R_\ep)$.

\bpr{DKP-lift} \cite{DKP1,Ha} For every $z \in Z(R_\ep)$, the hamiltonian derivation $y \mt \{ z, y \}$ of the Poisson algebra $( Z(R_\ep), \{. , .\})$ has a lift to an algebra derivation 
of $R_\ep$ given by 
\[
\partial_x ( \sig(\wt{y})) := \sig \left( \frac{x\wt{y} -\wt{y}x}{q-\ep} \right), \quad x \in \sig^{-1}(z), \wt{y} \in R.  
\]
\epr
Note that $\sig(x) \in Z(R_\ep)$ implies that $\sig (x\wt{y} - \wt{y}x) =0$, so $x\wt{y} - \wt{y}x \in (q-\ep) R$. 
The lifts coming from different elements $x, x' \in \sig^{-1}(z)$ differ by the inner derivation of $R_\ep$ corresponding to $\sig(  (x-x')/(q-\ep))$.

\bth{discr-symp} Let $R$ be a $\KK[q^{\pm1}]$-algebra for a field $\KK$ of characteristic 0 and 
$\ep \in \KK^\times$ be such that $q - \ep \in R$ is regular.
Assume that $R_\ep := R/(q-\ep)R$ is a free module of finite rank over a Poisson
subalgebra $C_\ep$ of its center. 

(i) Then $d(R_\ep/C_\ep)$ is a Poisson normal element of $(C_\ep, \{.,.\})$.

(ii) Assume, in addition, that $C_\ep$ is a unique factorization domain as a commutative algebra or a noetherian Poisson unique factorization domain. 
Then, $d(R_\ep/C_\ep)=0$ or
\[
d(R_\ep/C_\ep) =_{C_\ep^\times} \prod_{i=1}^m p_i
\]
for some (not necessarily distinct) Poisson prime elements $p_1, \ldots, p_m \in C_\ep$.
\eth
As usual, a product of 0 primes is considered to be 1.
Let $(A, \{.,.\})$ be a Poisson algebra and $u \in A^\times$. Then $a \in A$ is Poisson normal if and only if $ua$ is Poisson normal. The  
discriminant $d(R_\ep/C_\ep)$ is defined up to a unit of $C_\ep$, but because of this property it does not matter which representative 
is considered in part (i) of the theorem.

If $R_\ep$ is an order in a central simple algebra, then the discriminant $d(R_\ep/C_\ep)$ is nonzero. Specializations of 
iterated skewpolynomial extensions fall in this class by \cite[Theorem 1.5]{DKP2}. In particular, this is true for 
the families of algebras considered in the next two sections. The nonvanishing of the discriminants of those algebras 
also follows from the fact that these algebras have filtrations whose associated graded algebras are quasipolynomial 
algebras, see \eqref{filt}; by \cite[Proposition 4.10]{CPWZ2} the leading terms of the discriminants are nonzero.
Generally, nonvanishing of discriminants for Cayley--Hamilton algebras 
follows from the description of the kernel of the trace form in \cite[Proposition 3.4 (2)]{DPRR}.

By \cite[Example 5.12]{J}, there are examples of Poisson structures on polynomial algebras that are not Poisson UFDs. In the opposite
direction, it is easy to construct Poisson UFDs that are not UFDs as commutative algebras. In other words the two classes of algebras 
in \thref{discr-symp} (ii) are not properly contained in each other.
\subsection{}
\label{3.2} The next result is an explicit version of the statement in \thref{discr-symp} (i).
\bpr{Poisson-brack} In the setting of \thref{discr-symp} (i), let $\YY := \{ y_1, \ldots, y_N \}$ be a $C_\ep$-basis of $R_\ep$.
For all $z \in C_\ep$ and $x \in \sig^{-1}(z)$, we have 
\[
\{ z,  d_N(\YY : \tr) \} =  2 \tr(b(x)) d_N(\YY : \tr)
\]
where $b(x) := (b_{ij}) \in M_N(C_\ep)$ is the matrix with entries given by 
\[
\partial_x(y_i) := \sum_{j} b_{ij} y_j.
\]
\epr
\begin{proof} Set $\de:= d_N(\YY : \tr)$. The proposition follows by combining Propositions \ref{pinvar} and \ref{pDKP-lift}:
\[
\{ z, \delta \} = \partial_x \delta = 2 \tr(b(x)) \delta.
\]
\end{proof}
Part (i) of \thref{discr-symp} follows from \prref{Poisson-brack}. The second part follows from the first and 
\prref{fact}.
%%%%%%%%%
\subsection{}
\label{3.3}
In the situations in which the problem for computing the discriminant $d(R_\ep/C_\ep)$ was posed,
$C_\ep$ differs only slightly from the full center $Z(R_\ep)$. The restriction of the Poisson structure $\{.,.\}$ 
to $C_\ep$ is very nontrivial because of the nature of the definition in \eqref{Poisson}. This causes the
collection of Poisson primes of $C_\ep$ to be a small subset of the set of all prime elements of $C_\ep$. \thref{discr-symp} 
places a strong restriction on the possible form of the discriminant $d(R_\ep/C_\ep)$. One can fully determine it using the following 4 methods 
and sets of existing results from Poisson geometry and algebra:

{\bf{(1)}} If the algebra $R_\ep$ is $\Zset^n$-graded and $C_\ep$ is a homogeneous subalgebra, then one 
can choose a homogeneous $C_\ep$-basis $\YY$ of $R_\ep$. Since, in this case, the trace map $\tr \colon R_\ep \to C_\ep$
will be homogeneous, $d_N(\YY : \tr)$ will be graded and 
\[
\deg d_N(\YY : \tr) = 2 \sum_{y \in \YY} \deg y.
\]
Furthermore, the grading assumption implies $C_\ep^\times = (C_\ep)_0^\times$, 
thus the class of associates for $d(R_\ep/C_\ep)$ will consist of homogeneous elements  
of the same degree. The primes in \thref{discr-symp} (ii) will need to be homogeneous and their degrees will satisfy
\[
\sum_{i=1}^m \deg p_i = \deg d(R_\ep/C_\ep) = 2 \sum_{y \in \YY} \deg y.
\]

{\bf{(2) (A)}} The symplectic foliations of the Poisson manifolds coming up in the theory of quantum groups 
are well understood: the Belavin--Drinfeld Poisson structures \cite{Y-Duke}, the varieties of Lagrangian 
subalgebras \cite{EL1,EL2}, the Poisson homogeneous spaces of non-standard Poisson structures on 
simple Lie groups \cite{LY}. In light of \reref{prime-norm-def} (iii), these facts can be translated into 
results for the Poisson primes of the corresponding coordinate rings. The results will be for the case when the base field is $\Cset$, 
but the algebras in the theory of quantum groups are defined over $\Qset[q^{\pm1}]$ and by base 
change one can convert the results to any base field of characteristic 0.

{\bf{(B)}} The Poisson primes of all algebras in the very large class of so called Poisson--CGL extensions 
are described in \cite{GY}.

Combining (A) and (B), gives a description of the Poisson primes needed for \thref{discr-symp} (ii) for broad classes 
of algebras.

{\bf{(3)}} If $C_\ep$ is a domain, \thref{discr-symp} (i) implies that $d_N(\YY : \tr)$ gives rise to a derivation 
$\partial_{discr}$ of $C_\ep$ such that 
\[
\{d_N(\YY : \tr), z \} = d_N(\YY : \tr) \partial_{discr} (z), \quad \forall z \in C_\ep.
\]
This derivation is explicitly given by \prref{Poisson-brack}.
Every Poisson prime $p \in C_\ep$ also gives rise to a derivation $\partial_p$ of $C_\ep$ such that 
\[
\{p, z \} = p \partial_p (z), \quad \forall z \in C_\ep.
\]
The primes in \thref{discr-symp} need to satisfy
\[
\sum_{i=1}^m \partial_{p_i} = \partial_{discr}.
\]

In fact, the procedure (3) can be also applied to the general situation when the conditions in \thref{discr-symp} (ii) are not satisfied. 
\prref{Poisson-brack} determines the Poisson brackets of $d_N(\YY : \tr)$ with all hamiltonians on $\Spec C_\ep$ from which one can determine 
the evolution of $d_N(\YY : \tr)$ under all hamiltonian flows on $\Spec C_\ep$.
%and 
%From a geometric point of view, this matches the evolution of the discriminant under the hamiltonian flows on $\Spec C_\ep$ to the 
%evolution of the Poisson primes on $C_\ep$.

{\bf{(4)}} Filtrations of the algebra $R$ can be used to obtain leading term results for $d(R_\ep/C_\ep)$, see \cite[Proposition 4.10]{CPWZ2}.
They put further restrictions on what Poisson primes can appear in the expansion in \thref{discr-symp} (ii) by comparing the leading terms
of the two sides. In concrete situations these filtrations are different from the gradings in (1).

In the next 2 sections we show how one can use these methods to compute explicitly the discriminants
of the specializations at roots of unity of the algebras of quantum matrices and more generally those of all quantum Schubert cell algebras.
%%%%%%%%%
\subsection{}
\label{3.4} Next, we illustrate how \thref{discr-symp} and the recipe (1-4) give easy derivations for previous results for discriminants 
that were obtained by more involved methods.

The 2 dimensional quantum Weyl algebra is the $\KK[q^{\pm 1}]$-algebra $A_q$ with generators $x_1, x_2$ and  
relation 
\[
x_1 x_2 = q x_2 x_1 + 1.
\]
Let $l \in \Zset$, $l >1$ and $\ep$ be a primitive $l$-th root of unity in $\KK$.
Consider the specialization $A_\ep := A_q / (q-\ep) A_q$ and the canonical projection $\sig \colon A_q \to A_\ep$. The center
$Z(A_\ep)$ is isomorphic to the polynomial algebra in $z_1:= \sig(x_1)^l$ and $z_2:= \sig(x_2)^l$ over $\KK$. 
$A_\ep$ is a free $Z(A_\ep)$-module with basis 
\[
\YY:= \{ \sig(x_1)^{k_1} \sig(x_2)^{k_2} \mid k_1, k_2 \in [0,l-1] \}. 
\]
\bth{qWeyl} \cite{CYZ} For any base field $\KK$ of characteristic 0, $l \in \Zset$, $l >1$ and a primitive $l$-th root of unity $\ep \in \KK$,
the discriminant of the specialization of the quantum Weyl algebra $A_\ep$ over its center is given by
\[
d(A_\ep/Z(A_\ep)) =_{\KK^\times} (z_1 z_2 - t)^{l(l-1)} \quad \mbox{where} \quad t := (1-\ep)^{-l}.
\]
\eth
\begin{proof}
The induced Poisson bracket \eqref{Poisson} on $Z(A_\ep)$ is
\begin{equation}
\label{Pst}
\{ z_1, z_2 \} = l^2 \ep^{-1} (z_1 z_2-t).
\end{equation}
This follows from the expansion
\[
x_1^l x_2^l - q^{l^2} x_2^l x_1^l = \sum_{i=0}^{l-1} t_i x_2^i x_1^i 
\quad \mbox{with} \quad t_0 = \prod_{m=1}^l(1+ q + \cdots + q^{m-1}), t_i \in \KK[q^{\pm 1}]
\]
and the fact that $\sig(t_i x_2^i x_1^i/(q-\ep)) =0$ for $0 < i < l$ since $\{z_1, z_2\} \in Z(A_\ep) = \KK[z_1, z_2]$.

It is sufficient to prove the statement of the theorem for $\KK = \Cset$ since the structure constants for the 
products of the elements of $\YY$ and $z_i$ belong to $\Qset(\ep)$.
When $\KK = \Cset$, the symplectic leaves of the Poisson structure \eqref{Pst} are 
\begin{itemize}
\item
The complement to the hyperplane $z_1 z_2 = t$ (2 dimensional leaf),
\item The points on the hyperplane $z_1 z_2 = t$ (0 dimensional leaves).
\end{itemize}
By \reref{prime-norm-def} (iii), the only Poisson prime of $C_\ep$ is $z_1 z_2 - t$. Indeed, the zero locus
of any other irreducible polynomial $f \in \Cset[z_1, z_2]$ will intersect the 2-dimensional symplectic leaf nontrivially.
Hence, such an $f$ cannot be Poisson prime. (Since the Poisson bracket is defined over $\Qset$, it follows from 
\cite[Theorem 3.4]{GY} that $z_1 z_2 - t$ is the only Poisson prime element of $C_\ep$ for 
any base field $\KK$ of characteristic 0 such that $\ep \in \KK$.)
Since $Z(A_\ep)$ is a polynomial algebra and thus a UFD, \thref{discr-symp} implies
\[
d(A_\ep/Z(A_\ep)) =_{\KK^\times} (z_1 z_2 - t)^m
\]
for some $m \in \Nset$, We are left with computing $m$ which we do by using the method (3) in \S \ref{3.3}.
Since $\{ z_1, (z_1 z_2 - t) \} = l^2 \ep^{-1} z_1 ( z_1 z_2 - t)$, 
\begin{equation}
\label{eq}
\{ z_1, d_{l^2} ( \YY : \tr) \} =m l^2 \ep^{-1} z_1 d_{l^2} ( \YY : \tr).
\end{equation}
For $r \in A_\ep = \oplus_{y \in \YY} \KK[z_1, z_2] y$ and $y \in \YY$, denote by 
$\coeff_{z_1, y}(r)$ the coefficient of $z_1 y$ in $r$. One easily checks that for $y =  \sig(x_1)^{k_1} \sig(x_2)^{k_2}$,
\[
\coeff_{z_1, y} \big( \partial_{x_1^l}(y) \big) = k_2 l \ep^{-1}. 
\]
Propositions \ref{pDKP-lift} and \ref{pPoisson-brack}, and eq. \eqref{eq} imply
\[
m l^2 \ep^{-1}= \frac{\{ z_1, d_{l^2} ( \YY : \tr) \}}{z_1 d_{l^2} ( \YY : \tr)} = 2 \sum_{y \in \YY} \coeff_{z_1, y} \big( \partial_{x_1^l}(y) \big) = 
2 l^2 \ep^{-1} \sum_{k_2=0}^{l-1} k_2 = l^3(l-1) \ep^{-1}.
\]
So, $m = l(l-1)$.
\end{proof}

\bre{Domain} In \cite{CPWZ1,CPWZ2,CYZ} the more general problem of computing discriminants of algebras over integral domains $A$ 
was considered. One can obtain extensions of Theorems \ref{tdiscr-symp} (ii), \ref{tqWeyl} and the results below for specializations 
of algebras $R$ over $A[q^{\pm 1}]$ for an integral domain $A$ as follows. First, apply the theorems to the algebras 
$R \otimes_A Q(A)$ over $Q(A)[q^{\pm 1}]$ where $Q(A)$ is the field of fractions of $A$; this would compute the 
discriminants $d_N(R_\ep \otimes_A Q(A), C_\ep \otimes_A Q(A))$. Then compute the leading term of $d(R_\ep/C_\ep)$
over $A$ using \cite[Proposition 4.10]{CPWZ2}, i.e., step (4) in \S \ref{3.3}, and convert the formula for 
$d_N(R_\ep \otimes_A Q(A), C_\ep \otimes_A Q(A))$ to one for $d(R_\ep/C_\ep)$ by clearing the denominators and 
introducing the necessary extra factor from $A$ in $d(R_\ep/C_\ep)$.
\ere
%%%%%%%%%%%%
\subsection{}
\label{3.5}
Next we prove two general results for the $n$-discriminant ideals of specializations of algebras, 
see \deref{discr} (ii) and \cite[p. 126]{Re} for background on this notion. 
These results do not assume any freeness conditions like the one in \thref{discr-symp}. 
\bth{last} 
Let $R$ be a $\KK[q^{\pm1}]$-algebra for an infinite field $\KK$ and 
$\ep \in \KK^\times$ be such that $q - \ep \in R$ is regular.
Assume that $C_\ep$ is a Poisson subalgebra of the center of $R_\ep := R/(q-\ep)R$
and that $R_\ep$ is equipped with a trace function $\tr : R_\ep \to C_\ep$ which commutes with 
all derivations $\partial$ of $R_\ep$ such that $\partial(C_\ep) \subseteq C_\ep$.

Then, for all positive integers $n$, the discriminant ideal $D_n(R_\ep/C_\ep)$ is a 
Poisson ideal of $C_\ep$. Furthermore, it has the property that 
$\partial (D_n(R_\ep/C_\ep)) \subseteq D_n(R_\ep/C_\ep)$ for all 
derivations $\partial$ of $R_\ep$ such that $\partial(C_\ep) \subseteq C_\ep$.
\eth

The first statement in the theorem follows from the second in view of \prref{DKP-lift}.
The second statement of the theorem follows from the next proposition.

\bpr{derI} Assume that $\tr \colon S \to C \subset Z(S)$ is a trace for an algebra $S$ over an infinite field $\KK$ which
commutes with all derivations $\partial$ of $S$ such that $\partial(C) \subseteq C$.
Then
$\partial(D_n(S/C)) \subseteq D_n(S/C)$ for all derivations $\partial$ of $S$ such that $\partial(C) \subseteq C$.
\epr
Given a positive integer $n$, define 
\[
\lcor .,. \rcor \colon S^n \times S^n \to C \quad \mbox{by} \; \; 
\lcor \XX, \YY \rcor := \det \big( [\tr(x_i y_j)]_{i,j=1}^n \big)
\]
for $\XX:= (x_1, \ldots, x_n)$, $\YY:= (y_1, \ldots, y_n) \in S^n$. This is 
obviously a symmetric form on $S^n$ which is $C$-polylinear in the sense that 
$\lcor (x_1, \ldots, c x_k, \ldots, x_n), \YY \rcor = c \lcor \XX, \YY \rcor$ for
all $c \in C$ and $k \in [1,n]$. For a derivation $\partial$ of $S$, define
\[
\ol{\partial}(\XX):= ( \partial (x_1), \ldots, \partial(x_n) ) 
\]
and
\[
\partial(\XX):= \sum_{k=1}^n (x_1, \ldots, \partial (x_k), \ldots, x_n).
\]
\noindent
{\em{Proof of \prref{derI}.}} For $p(t) \in S[t]$, denote by $\coeff_{t^i} p(t) \in S$ 
the coefficient of $t^i$ in $p(t)$. Using several times the differentiation property of $\partial$ and the assumption that $\partial$ commutes with $\tr$, gives
\[
\partial( d_n(\XX : \tr) ) = 2 \lcor \XX, \partial(\XX) \rcor = \coeff_t (d_n(\XX + t \ol{\partial}(\XX) : \tr))
\]
for all $\XX \in S^n$. The proposition follows from the fact that $d_n(\XX + t \ol{\partial}(\XX) : \tr) \in D_n(S/C),$ $\forall t \in \KK$
and the assumption that $\KK$ is infinite.
\qed

\thref{last} and \prref{derI} have natural bilinear analogs. 
Let $S$ be an algebra with trace $\tr \colon S \to C$ where $C$ is a subalgebra of $Z(S)$. 
Following \cite[Definition 1.2 (2)]{CPWZ2},
define the $n$-th {\em{modified discriminant ideal}} $MD_n(S/C)$ to be the ideal of $C$, generated by
\[
\lcor \XX, \YY \rcor \quad \mbox{for all} \; \;  \XX, \YY \in S^n. 
\] 
Thus, $D_n(S/C) \subseteq MD_n(S/C)$. If $S$ is a free rank $N$ module over $C$ with a basis $\XX \in S^N$,
then 
\[
MD_N(S/C) = D_N(S/C) = (d_N(\XX : \tr))
\]
by an argument similar to the identity \eqref{discr-bas}. We refer the 
reader to \cite[Sect. 1]{CPWZ2} for other properties of modified discriminant ideals.

\bth{bilin} Assume that $\KK$ is an infinite field and $n$ is a positive integer.

(i) Let $\tr \colon S \to C \subset Z(S)$ be a trace for a $\KK$-algebra $S$ which commutes with 
all derivations $\partial$ of $S$ such that $\partial(C) \subseteq C$.
Then $\partial(MD_n(S/C)) \subseteq MD_n(S/C)$ 
for all derivations $\partial$ of $S$ such that $\partial(C) \subseteq C$.

(ii) In the setting of \thref{last}, for all $\ep \in \KK^\times$, 
the modified discriminant ideal $MD_n(R_\ep/C_\ep)$ is a Poisson ideal of $C_\ep$ with respect to the 
induced Poisson structure on $C_\ep$. Moreover, $\partial (MD_n(R_\ep/C_\ep)) \subseteq MD_n(R_\ep/C_\ep)$ for all 
derivations $\partial$ of $R_\ep$ such that $\partial(C_\ep) \subseteq C_\ep$.
\eth
\thref{bilin} (i) is proved analogously to \prref{derI} using the identity
\[
\partial \lcor \XX, \YY \rcor = \lcor \XX, \partial(\YY) \rcor + 
\lcor \partial(\XX), \YY \rcor = \coeff_{t} \lcor \XX + t \ol{\partial}(\XX), \YY + t \ol{\partial}(\YY) \rcor
\]
for all $\XX, \YY \in S^n$, obtained by applying the differentiation property of $\partial$ and 
the assumption that $\partial$ commutes with $\tr$. The second 
part of the theorem follows from the first. 
\subsection{}
\label{3.6}
We finish the section with a generalization of the results in \S \ref{3.1} and \S \ref{3.5} to the framework of Poisson orders introduced by Brown and Gordon.

\bde{Pois-order} \cite{BrGo}
Assume that $S$ is an affine algebra over a field $\KK$ of characteristic 0 which is a finite module over a central subalgebra 
$C$. The algebra $S$ is called a Poisson $C$-order if there is a $\KK$-linear map $\De \colon C \to \Der_\KK(S)$ such that
\begin{enumerate}
\item $C$ is stable under $\De_z$ for all $z \in C$ and
\item the induced bracket $\{.,.\}$ on $C$ given by $\{z_1, z_2 \}:= \De_{z_1}(z_2)$ turns $C$ into a Poisson algebra.
\end{enumerate} 
\ede
\prref{DKP-lift} implies that, in the setting of the proposition, $R_\ep$ is a Poisson $Z(R_\ep)$-order. The map $\De$ is given as follows. 
Choose a $\KK$-linear map $\om \colon Z(R_\ep) \to R$ such that $\sig \om = \id_{Z(R_\ep)}$ and set $\De_z = \partial_{\om(z)}$.

\bth{P-order} Let $S$ be a $\KK$-algebra over a field $\KK$ (of characteristic 0) which is 
a Poisson $C$-order. Let $\tr \colon S \to C$ be a trace map that commutes with all derivations 
of $S$ that preserve $C$. 

(i) If $S$ is a free $C$-module (of finite rank), then $d(S/C)$ is a Poisson normal element of $C$. If, in addition, $C$ is a unique factorization domain 
as a commutative algebra or a noetherian Poisson unique factorization domain, then either $d(S/C)=0$ or
\[
d(S/C) =_{C^\times} \prod_{i=1}^m p_i
\]
for some (not necessarily distinct) Poisson prime elements $p_1, \ldots, p_m \in C$.

(ii) For all positive integers $n$, the discriminant and modified discriminant ideals $D_n(S/C)$ 
and $MD_n(S/C)$ are Poisson ideals of $C$. Furthermore, $\partial (D_n(S/C)) \subseteq D_n(S/C)$ 
and $\partial (MD_n(S/C)) \subseteq MD_n(S/C)$ for all 
derivations $\partial$ of $S$ such that $\partial(C) \subseteq C$.
\eth
The theorem follows from Propositions \ref{pinvar}, \ref{pfact} and \ref{pderI} and \thref{bilin} (i).
We leave the details to the reader. 
%%%%%%%%%%%%%%%%
\sectionnew{Discriminants of quantum matrices}
\label{q-matr}
In this section we derive a formula for the discriminants of the specializations at odd roots of unity of the 
algebras of square quantum matrices. The results do not require any Lie theoretic background and are proved independently 
from the results of De Concini, Kac and Procesi \cite{DKP1,DKP1b} on the specializations of quantized universal enveloping algebras 
at roots of unity. The case of rectangular quantum matrices can be obtained along the same lines, but is somewhat more technical. 
Because of this, it appears in \S \ref{5.5} and is obtained by applying the Lie theoretic results from the next section.
\subsection{}
\label{4.1}
Throughout this section, we fix a field $\KK$ of characteristic 0 and an odd positive integer $l>2$ such that 
$\KK$ contains a primitive $l$-th root of unity $\ep$.

Let $n \in \Zset$, $n>1$. The algebra of square quantum matrices $R_q[M_n]$ is the $\KK[q^{\pm 1}]$-algebra with generators 
$x_{ij}$, $i, j \in [1,n]$ and relations 
\begin{align}
\nn
x_{ij} x_{kj} &= q x_{kj} x_{ij}, \quad \mbox{for} \; i < k, \\
\label{qmatr}
x_{ij} x_{ir} &= q x_{ir} x_{ij}, \quad \mbox{for} \; j < r, \\ 
\nn
x_{ij} x_{kr} &= x_{kr} x_{ij}, \quad 
\mbox{for} \; i < k, j > r,\\
\nn
x_{ij} x_{kr} - x_{kr} x_{ij} &= (q-q^{-1}) x_{ir} x_{kj}, \quad 
\mbox{for} \; i < k, j<r.
\end{align}
Denote the specialization $R_\ep[M_n] :=  R_q[M_n] / (q - \ep) R_q[M_n]$ and the natural projection
$\sig \colon R_q[M_n] \to R_\ep[M_n]$. 
One easily verifies that the following elements are central
\[
z_{ij}:= \sig(x_{ij})^l \in Z(R_\ep[M_n]).
\]
Denote by $C_\ep[M_n]$ the subalgebra of $Z(R_\ep[M_n])$ generated by them. It is isomorphic to a polynomial 
algebra: 
\begin{equation}
\label{centr}
C_\ep[M_n] \cong \KK[z_{ij}, 1 \leq i,j \leq n].
\end{equation}
For $k \in [1,n]$, denote the minors 
\[
\De_k := \De_{[n-k+1, n]; [1, k]}, \quad 
\ol{\De}_j:= \De_{[1, j] ; [n-j+1, n]} \in C_\ep[M_n].
\] 
As in the introduction, $\De_{I;J}$ denotes the minor in $z_{ij}$ for the rows in the set $I \subseteq [1,n]$ and the columns in $J\subseteq [1,n]$. 
Note that $\De_n= \ol{\De}_n$.

The algebra $R_\ep[M_n]$ is a free module over $C_\ep[M_n]$ with basis
\begin{equation}
\label{q-matr-basis}
\YY:=\{ \sig(x_{11})^{k_{11}} \ldots \sig(x_{nn})^{k_{nn}} \mid 0 \leq k_{11}, \ldots, k_{nn} \leq l-1 \}
\end{equation}
where the (noncommuting) $\sig(x_{ij})$-elements are listed down the first column, second, etc., to the $n$-the column. The facts 
in \eqref{centr} and \eqref{q-matr-basis} hold because $\{ x_{11}^{k_{11}} \ldots x_{nn}^{k_{nn}} \mid k_{11}, \ldots, k_{nn} \in \Nset \}$ is
a $\KK[q^{\pm 1}]$-basis of $R_q[M_n]$.

\bth{q-matr-discr} Let $\KK$ be a field of characteristic 0, $l >2$ an odd integer 
and $\ep \in \KK$ a primitive $l$-th root of unity. Then, 
\[
d(R_\ep[M_n]/C_\ep[M_n]) =_{\KK^\times} \prod_{k=1}^n \De_k^L
\prod_{j=1}^{n-1} \ol{\De}_j^{\, L}
\] 
where $L:= l^{n^2-1}(l-1)$. 
\eth
The theorem is proved in \S \ref{4.3} and \S \ref{4.2} contains the needed results for the related Poisson structures.  
%%%%%%%%%%%%%
\subsection{}
\label{4.2}
\ble{ind-Poisson} For an odd integer $l>2$, $C_\ep[M_n]$ is a Poisson subalgebra of $Z(R_\ep[M_n])$ with respect to the Poisson 
structure \eqref{Poisson} and the corresponding Poisson bracket is given by 
\begin{equation}
\label{P-brack}
\{ z_{ij}, z_{km} \} := l^2 \ep^{-1} (\sign(k - i )+\sign(m - j )) z_{im} z_{kj}.
\end{equation}
\ele
The lemma can be derived from \cite[Theorem 7.6]{DKP1}. We sketch a direct proof below that shows 
at what stages the condition of $l$ being odd is needed.
\begin{proof} For all $1 \leq i<k \leq n$, $1 \leq j < m \leq n$, we have the algebra embeddings $R_q[M_2] \hra R_q[M_n]$ given by 
\begin{equation}
\label{embed-qm}
\begin{bmatrix}
x_{11} & x_{12} \\
x_{21} & x_{22} 
\end{bmatrix}
\mt
\begin{bmatrix}
x_{ij} & x_{im} \\
x_{kj} & x_{km} 
\end{bmatrix}.
\end{equation}
Thus, it suffices to check the statement for $n=2$. In that case the only identity that is not immediate is
\[
\{ z_{11}, z_{22} \} = 2 l^2 \ep^{-1} z_{12} z_{21}.
\]
To obtain this, we first compute the derivation
\begin{align*}
\partial_{x_{11}^l} (\sig(x_{22})) &= \sig \left( \frac{[x_{11}^l, x_{22}]}{q-\ep}\right) = \sig \left( \sum_{i=0}^{l-1} \frac{x_{11}^i [x_{11}, x_{22}] x_{11}^{l-i-1}}{q-\ep}\right)
\\
&= \sig \left( \frac{q^{2l} -1}{q(q - \ep)} x_{12} x_{21} x_{11}^{l-1} \right) = 2 l \ep^{-2} \sig(x_{12}) \sig(x_{21}) \sig(x_{11})^{l-1}.
\end{align*}
This is valid for all integers $l>2$. Applying \prref{DKP-lift} gives
\begin{equation}
\label{zbrack}
\{ z_{11}, z_{22} \} = \partial_{x_{11}^l} (\sig(x_{22})^l)  = 2 l \ep^{-2} \sum_{i=0}^{l-1} \sig(x_{22})^i \sig(x_{12}) \sig(x_{21}) \sig(x_{11})^{l-1} \sig(x_{22})^{l-i-1}.
\end{equation}
From the commutation relations in $R_q[M_n]$ we obtain that the right hand side can be uniquely represented in the 
form 
\[
\{ z_{11}, z_{22} \}  = \sum_{i=1}^{l} c_i \sig(x_{12})^i \sig(x_{21})^i \sig(x_{22})^{l-i}  \sig(x_{11})^{l-i}
\]
for some $c_i \in \KK$.

Since $\{z_{11}, z_{22} \} \in Z(R_\ep[M_2])$, it commutes with $\sig(x_{22})$. This gives the recursion relation 
\[
\ep^{2i-1} (\ep^{2 (l-i+1)} -1) c_{i-1} = (\ep^{2i} -1)c_i 
\]
from which we get $c_{l-1} = c_{l-2} = \ldots = c_1 = 0$. Here we use in an essential way the condition that $l$ is odd 
which gives that $\ep^{2i} - 1 \neq 0$ for $i \in [1,l-1]$.
The coefficient $c_l$ is computed directly from \eqref{zbrack}:
\[
c_l = 2l \ep^{-1} \prod_{i=1}^{l-1} (\ep^{2i}-1) = 2 l^2 \ep^{-1}.
\]

Here again we used that $l$ is odd, so $\ep^2$ is a primitive $l$-th root of unity and thus $\prod_{i=1}^{l-1} (t - \ep^{2i}) = \sum_{i=0}^{l-1} t^i$.
\end{proof}
\bre{evenl} In the case of even $l$, $C_\ep[M_n]$ is not a Poisson subalgebra of $Z(R_\ep[M_n])$. For example, for $l=4$, $\ep = \pm i$,
\[
\{z_{11}, z_{22} \} = -32 \ep \sig(x_{22})^2 \sig(x_{12})^2 \sig(x_{21})^{2}  \sig(x_{11})^{2}.
\]
\ere

The algebra $R_q[M_n]$ is $\Zset^{2n}$-graded by setting $\deg x_{ij} := e_i + e_{n+j}$ where $\{e_1, \ldots e_{2n} \}$ is 
the standard basis of $\Zset^{2n}$. This induces $\Zset^{2n}$-grading on $R_\ep[M_n] = R_q[M_n] /(q-\ep)R_q[M_n]$. The subalgebra $C_\ep[M_n]$ is homogeneous
($\deg z_{ij} := l (e_i + e_{n+j})$) and the trace map $\tr \colon R_\ep[M_n] \to C_\ep[M_n]$ is graded. Denote the symmetric bilinear form 
\[
\lcor .,. \rcor \colon \Zset^{2n} \times \Zset^{2n} \to \Zset, \quad \mbox{given by} \; \; 
\lcor e_i, e_j \rcor = \delta_{ij}.
\]
For $i \leq j$, denote
\[
e_{[i,j]} = e_i + \cdots + e_j.
\]

\bpr{q-matr-hpr} \cite{GY} The homogeneous prime elements of $C_\ep[M_n]$ with respect to the Poisson bracket
\eqref{P-brack} are 
\[
\De_1, \ldots, \De_n, \ol{\De}_1, \ldots , \ol{\De}_{n-1}.
\]
They satisfy
\[
\{ \De_k, z_{ij} \} = l^2 \ep^{-1} ( \lcor e_{[n-k+1,n]}, e_i \rcor - \lcor e_{[n+1, n+k]}, e_{n+j} \rcor ) \De_k z_{ij} 
\]
and 
\[
\{ \ol{\De}_k, z_{ij} \} = - l^2 \ep^{-1} ( \lcor e_{[1,k]}, e_i \rcor - \lcor e_{[2n-k+1, 2n]}, e_{n+j} \rcor) \ol{\De}_k z_{ij} 
\]
for all $i,j \in [1,n]$.
\epr
In the case $\KK = \Cset$, \prref{q-matr-hpr} has a very illuminating proof using Poisson geometry which we give below. 
This is sufficient to deduce the case of arbitrary base field in \thref{q-matr-discr} because the structure constants for the 
products of the elements in $\YY$ and $z_{ij}$ are in $\Qset(\ep)$; that is, proving \thref{q-matr-discr} over $\Cset$, will imply the 
statement over any base field of characteristic 0 such that $\ep \in \KK$.
\begin{proof} Let $\Cset = \KK$.
We identify $C_\ep[M_n] \cong \Cset[z_{ij}, 1 \leq i,j \leq n ]$ with the coordinate ring of the affine space $M_n$ of $n \times n$ complex 
matrices. The action of $H := (\Cset^\times)^{2n}$ on $M_n$,  given by
\begin{equation}
\label{Hmatr}
(t_1, \ldots, t_{2n}) \cdot [z_{ij}]_{i,j=1}^n := [t_i t_{n+j} z_{ij}]_{i,j=1}^n,
\end{equation}
preserves the Poisson structure \eqref{P-brack}. An element $f \in C_\ep[M_n]$ which is not a power of another element
is a homogeneous Poisson prime if and only if its zero 
locus is an irreducible subvariety of $M_n$ which is a union of $H$-orbits of symplectic leaves. By \cite[Theorem 0.4 (a)-(b)]{BGY} the $H$-orbits of symplectic 
leaves of $M_n$ are of the following 3 types:
 \begin{itemize}
 \item A Zariski dense $H$-orbit of leaves $\LL$.
 \item $2n-1$ $H$-orbits of leaves $\LL_1, \ldots, \LL_{2n-1}$ of codimension $1$. They are given by
 \[
 \LL_i =  \Big\{ g \in M_n \mid 
 \begin{bmatrix}
 I_n & w_n^\circ g \\
 0 & I_n
 \end{bmatrix} \in B_- s_i B_+ c^n \Big\}.
 \]  
 and are locally closed, irreducible. Here $B_\pm \subset {\mathrm{GL}}_{2n}(\Cset)$ are the standard Borel subgroups 
 of upper and lower triangular matrices, $s_i$ are the simple reflections of $S_{2n}$ realized as permutation matrices, 
 $w_n^\circ \in S_n \subset S_{2n}$ is the longest element of $S_n$, $c = (1 2 \ldots (2n))$ is the standard Coxeter element of $S_{2n}$.   
 \item Finitely many $H$-orbits of leaves of codimension $>1$.
 \end{itemize}
By \cite[Theorem 5.3]{Y-P} the zero loci of $\De_1, \ldots, \De_n, \ol{\De}_{n-1}, \ldots, \ol{\De}_1$ are precisely the Zariski closures of $\LL_1, \ldots, \LL_{2n-1}$. 
Since $\LL_i$ are irreducible and are $H$-orbits of symplectic leaves, the elements $\De_i, \ol{\De}_j$ are Poisson primes. The Poisson brackets 
in \prref{q-matr-hpr} are well known, see \cite[Lemma 3.2]{GSV}. 

Assume that  $f \in C_\ep[M_n]$ is any other homogeneous Poisson prime. Its zero locus $\VV(f)$ cannot intersect $\LL$, since $\LL$ 
is a single $H$-orbit of leaves and this would imply that $\VV(f) \supset \ol{\LL} = M_n$. (Here and in \S \ref{5.2}--\ref{5.3}, $\ol{X}$ denotes the Zariski 
closure of a locally closed subset $X$.)
Because the third type of $H$-orbits of symplectic leaves of $M_n$ are
of codimension $>1$, there exists $i \in[1,2n-1]$ such that $\VV(f) \cap \LL_i \neq \varnothing$. Using one more time that $\LL_i$ 
is a single $H$-orbit of leaves implies 
$\VV(f) \supseteq \ol{\LL_i}$. Since $f$ is prime, this is only possible if $\VV(f) = \ol{\LL_i}$. Thus, $f$ is one of the Poisson primes $\De_i,\ol{\De}_j$.
\end{proof} 
%%%%%%%%%%%%
\subsection{}
\label{4.3} The idea of the proof of \thref{q-matr-discr} is to apply \thref{discr-symp} and use the methods (1-3) in \S 3.3: We use the 
degree of $d(R_\ep[M_n]/C_\ep[M_n])$ and its Poisson brackets with the elements  of $C_\ep[M_n]$
which are obtained from \prref{Poisson-brack}. Denote for brevity
\[
\delta := d_{l^{n^2}}(\YY : \tr)
\]
for the $C_\ep[M_n]$-basis $\YY$ of $R_\ep[M_n]$ given by \eqref{q-matr-basis}. 
$C_\ep[M_n]$ is a polynomial ring and thus a UFD.
\thref{discr-symp} and \prref{q-matr-hpr} imply
\begin{equation}
\label{mij}
\delta =  t \prod_{i=1}^n \De_i^{m_i} \prod_{j=1}^{n-1} \ol{\De}_j^{\, \ol{m}_j}
\end{equation}
for some $t \in C_\ep[M_n]^\times = \KK^\times$ and $m_i, \ol{m}_j \in \Nset$.
Since $\tr \colon R_\ep[M_n] \to C_\ep[M_n]$ is $\Zset^{2n}$-graded, 
\[
\deg \delta = 2 \sum_{k_{11}, \ldots, k_{nn}=0}^{l-1} \deg x_{11}^{k_{11}} \ldots x_{nn}^{k_{nn}} =
n l^{n^2}(l-1) e_{[1,2n]}.
\]
Taking into account that 
\[
\deg \De_k = l (e_{[n-k+1,n]} + e_{[n+1, n+k]}), \quad \deg \ol{\De}_k = l( e_{[1,k]} +  e_{[2n-k+1, 2n]})
\]
and comparing the degrees of the two sides of \eqref{mij}, gives
\begin{equation}
\label{m-1eq}
m_1 + \cdots + m_n = nl^{n^2-1}(l-1), \quad \ol{m}_i = m_{n-i}, \; \; i \in [1,n-1].
\end{equation}

Since $C_\ep[M_n] \cong \KK[z_{ij}, 1 \leq i,j \leq n]$ and $R_\ep[M_n] = \oplus_{y \in \YY} C_\ep[M_n] y$, every 
element $r \in R_\ep[M_n]$ can be uniquely represented in the form 
\[
r = \sum_{y \in \YY} \sum_{\mu} t_{\mu, y} \mu y
\]
for some $t_{\mu,y} \in \KK$
where the second sum runs over all monomials $\mu$ of $\KK[z_{ij}, 1 \leq i,j \leq n]$. For a monomial $\mu$ in $z_{ij}$ and 
$y \in \YY$, denote
\[
\coeff_{\mu, y}(r) := t_{\mu, y}.
\]

\ble{q-matr-aux} For all $y = \sig(x_{11})^{k_{11}} \ldots \sig(x_{nn})^{k_{nn} } \in \YY$, 
\begin{equation}
\label{coeff}
\coeff_{z_{ij}, y} \big( \del_{x_{ij}^l} (y) \big) =  
\Big( \sum_{a=1}^n \sign(a-j) k_{ia} + \sum_{a=1}^n \sign(a-i) k_{aj} \Big) l \ep^{-1}.
\end{equation}
In particular, 
\[
\sum_{y \in \YY} \coeff_{z_{ij}, y} \big( \del_{x_{ij}^l} (y) \big) = (n-i-j+1)(l-1)l^{n^2+1} \ep^{-1}.
\]
\ele
\begin{proof} For ${\bf{k}} = (k_{11}, \ldots, k_{nn}) \in \Nset^{n^2}$, denote 
\begin{equation}
\label{elem-R}
x^{\bf{k}}:= x_{11}^{k_{11}} \ldots x_{nn}^{k_{nn} }
\end{equation}
where the product is ordered down columns 1, 2, etc., of the entries of an $n\times n$ matrix. These elements 
form a $\KK[q^{\pm 1}]$-basis of $R_q[M_n]$.
Eq. \eqref{coeff} in the lemma follows from the definition of the derivations $\partial_x$ in \prref{DKP-lift} and the following 
commutation relations which are obtained from the defining relations for $R_q[M_n]$:
\begin{align*}
& x_{ij}^l x^{\bf{k}} = q^{- \sum_{a=1}^{j-1} k_{ia} l - \sum_{a=1}^{i-1} k_{aj}l }  x^{\bf{k}'} + \; \mbox{lower order terms},
\\
& x^{\bf{k}} x_{ij}^l = q^{- \sum_{a=j+1}^n k_{ia} l - \sum_{a=i+1}^n k_{aj} l } x^{\bf{k}'} +  \; \mbox{lower order terms}
\end{align*}
where ${\bf{k}}' = (k'_{11}, \ldots, k'_{nn})$ with $k'_{ab} := k_{ab} + l \delta_{ai} \delta_{bj}$. The elements \eqref{elem-R}
are ordered using the (linear) lexicographic order on $\Zset^{n^2}$ from right to left, meaning that 
$(k_1, \ldots, k_{n^2}) \prec (m_1, \ldots, m_{n^2})$ if $k_{n^2} = m_{n^2}$, \ldots, $k_{i+1} = m_{i+1}$ and
$k_i < m_i$. The second equation in the lemma follows directly from the first by summation.
\end{proof}
To determine the integers $m_i$, $\ol{m}_j$, we use the Poisson brackets of $\de$ with the generators of $C_\ep[M_n]$, 
keeping in  mind  \thref{discr-symp} (i). \prref{q-matr-hpr} and eq. \eqref{mij}
imply that $\{ z_{ij}, \de \}/ (z_{ij} \de) \in \KK$. Combining this with \prref{invar}, gives
\[
2 \sum_{y \in \YY} \coeff_{z_{ij}, y} \big( \partial_{x_{ij}^l} (y) \big) z_{ij} \de = \{ z_{ij}, \de \} = M_{ij} z_{ij} \de
\]
where
\[
\frac{M_{ij}}{l^2 \ep^{-1}} = - \sum_{k=n-i+1}^n m_k + \sum_{k=i}^{n-1} \ol{m}_k + \sum_{k=j}^n m_k - \sum_{k=n-j+1}^{n-1} \ol{m}_k.
\]
Since $\ol{m}_i = m_{n-i}$, $M_{ij}$ simplifies to
\[
\frac{M_{ij}}{l^2 \ep^{-1}}
= \sum_{k=1}^{n-i} m_k - \sum_{k=n-i+1}^n m_k -  \sum_{k=1}^{j-1} m_k + \sum_{k=j}^n m_k.
\]
Thus, for $j \in [2,n]$,
\begin{align*}
2 m_j l^2 \ep^{-1} &= M_{i,j-1} - M_{ij} = 2 \sum_{y \in \YY} \coeff_{z_{i,j-1}, y} \big(\partial_{x_{i,j-1}^l} y \big) - 
2 \sum_{y \in \YY} \coeff_{z_{ij}, y} \big( \partial_{x_{ij}^l} y \big)
\\
& =  2 (l-1) l^{n^2+1} \ep^{-1}
\end{align*}
where at the end we used the second part of \leref{q-matr-aux}. Finally, $m_1 = (l-1) l^{n^2-1}$ from \eqref{m-1eq}.
This proves \thref{q-matr-discr}.
\sectionnew{Discriminants of quantum Schubert cell algebras}
\label{disrc-qSchubert}
In this section we prove an explicit formula for the discriminants of the specializations at roots of unity of the quantum Schubert cell algebras 
for all simple Lie algebras $\g$ and Weyl group elements $w \in W$.

For two subgroups $B_1$ and $B_2$ of a group $G$, we will denote by $g \cdot B_2$ the elements of $G/B_2$, by $B_1 g \cdot B_2$ 
the $B_1$-orbit of $g \cdot B_2 \in G/B_2$, and by $B_1 g B_2$ the corresponding double coset in $G$.
\subsection{}
\label{5.1} Let $[c_{ij}] \in M_r(\Zset)$ be a Cartan matrix of finite type and $\UU_q(\g)$ be the corresponding quantized universal enveloping algebra defined 
over $\KK(q)$ where $\KK$ is a field of characteristic 0. We will follow the notation of \cite{Ja} except for denoting
the Chevalley generators of $\UU_q(\g)$ by $E_i, F_i , K_i^{\pm 1}$, $i \in [1,r]$ (in \cite{Ja} they were indexed by the simple 
roots of $\g$). 

Let $W$ be the Weyl group of $[c_{ij}]$, $\{\al_1, \ldots, \al_r\}$ the set of simple roots, and $\{s_1, \ldots, s_r\} \subset W$ the corresponding set 
of simple reflections. Denote by $\lcor ., . \rcor$ the $W$-invariant bilinear form on $\bigoplus_j \Qset \al_i$ normalized by $\|\al_i\|^2 =2$ for short roots $\al_i$.
Set $q_i := q^{\| \al_i\|/2}$.

Given a Weyl group element $w$ and a reduced expression
\[
w = s_{i_1} \ldots s_{i_N},
\]
consider the roots $\be_k = s_{i_1} \ldots s_{i_{k-1}}(\al_{i_k}),$ $ k \in [1,N]$. They are precisely the roots of the nilpotent 
Lie algebra $\n_+ \cap w(\n_-)$ where $\n_\pm$ are the nilradicals of a pair of opposite Borel subalgebras of $\g$.  The quantum Schubert cell 
algebras $\UU^-[w]$ are the $\KK[q^{\pm1}]$-subalgebras of  $\UU_q(\g)$ generated by the quantum root vectors
\begin{equation}
\label{root-vect}
F_{\be_j} := T_{i_1} \ldots T_{i_{j-1}}(F_{i_j}), \quad j \in [1,N]
\end{equation}
where $T_i$ refers to the action \cite{Ja,L} of the braid group of $W$ on 
$\UU_q(\g)$. The algebras $\UU^-[w]$ do not depend on the choice of a reduced 
expression of $w$, \cite{L,DKP2}. Their generators satisfy the Levendorskii--Soibelman 
straightening relations: for $1 \leq j < m \leq N$, 
\begin{equation}
\label{rel}
F_{\be_m} F_{\be_j} - q^{ - \lcor \be_m, \be_j \rcor } F_{\be_j} F_{\be_m}  
= \sum_{k_{j+1}, \ldots, k_{m-1} \in \Nset} t_{k_{j+1}, \ldots, k_{m-1}} F_{\be_{j+1}}^{k_{j+1}} \ldots F_{\be_{m-1}}^{k_{m-1}} 
\end{equation}
for some $t_{k_{j+1}, \ldots, k_{m-1}} \in \Qset[q^{\pm1}]$. As a consequence, $\UU^-[w]$ has the PBW basis
\begin{equation}
\label{PBW}
\big\{F_{\be_1}^{k_1} \ldots F_{\be_N}^{k_N} \mid k_1, \ldots, k_N \in \Nset \big\}.
\end{equation}

\bre{OverQ} The algebras $\UU^-[w]$ can be defined as the algebras with generators $F_{\be_1}, \ldots, F_{\be_N}$ and 
relations \eqref{rel}. In particular, they are defined over $\Qset[q^{\pm1}]$ and their specializations at a root of unity $\ep$ 
are defined over $\Qset(\ep)$. All formulas for discriminants proved for one field of characteristic $0$ are valid for any other 
field of characteristic 0 by a direct base change.
\ere
%Following \cite{DKP1,DKP2}, define that, 
%\medskip
%
%{\em{$l >1$ is a good integer if it is odd and, in type $G_2$, $l > 3$.}} 
%\medskip

Let $\ep \in \KK$ be a primitive $l$-th root of unity.  
Denote the specialization $\UU^-_\ep[w]:= \UU^-[w]/(q-\ep)\UU^-[w]$ and the 
canonical projection $\sig \colon \UU^-[w] \to \UU^-_\ep[w]$. Set 
\begin{equation}
\label{zz}
z_{\be_j}:= (\ep^{ \| \al_{i_j} \| /2} - \ep^{-\|\al_{i_j} \| /2})^l \sig(F_{\be_j})^l \in \UU^-_\ep[w], \quad j \in [1,N].
\end{equation}
Denote by $C^-_\ep[w]$ the $\KK$-subalgebra of $\UU^-_\ep[w]$ generated by $z_{\be_j}$, $j \in [1,N]$.
 
\bth{cent} \cite{DKP1} For all integers $l>1$, $C^-_\ep[w]$ is a subalgebra of $Z( \UU^-_\ep[w])$. It is isomorphic to the polynomial 
algebra in the generators $z_{\be_j}$, $j \in [1,N]$ and is independent of the choice of reduced expression of $w$.
\eth
The last part was stated in \cite[Proposition 3.3]{DKP1} for the longest element of $W$; the proof works for all $w \in W$. The algebra $\UU^-_\ep[w]$
is a free $C^-_\ep[w]$-module with basis
\begin{equation}
\label{Uwbasis}
\YY:= \{\sig(F_{\be_1})^{k_1} \ldots \sig(F_{\be_N})^{k_N} \mid k_1, \ldots, k_N \in [0,l-1] \}.
\end{equation}
This and the second part of \thref{cent} follow from the PBW basis \eqref{PBW}. 

Denote by $G$ the split, connected, simply connected algebraic $\KK$-group with Lie algebra $\g$. 
Let $B_\pm$ be a pair of opposite Borel subgroups of $G$
and $U_\pm$ be their unipotent radicals. 
Let $\{e_i, f_i\}$ be a set of Chevalley generators of $\g$ that generate $\Lie(U_\pm)$. 
Denote by $\dot{s}_i$ the representatives of $s_i$ in the normalizer of the maximal torus $H: = B_+ \cap B_-$ of $G$ given by 
\[
\dot{s}_i := \exp(f_i) \exp (-e_i) \exp (f_i). 
\]
They are extended (in a unique way) to Tits' representatives of the elements $u \in W$ in $N_G(H)$ by setting $\dot{v} := \dot {u} \dot{s}_i$ if 
$v = u s_i$ and $\ell(v) = \ell(u) + 1$ where $\ell \colon W \to \Nset$ is the length function.  
Given a positive root $\be = u(\al_i)$ of $\g$, for some $u \in W$ and a simple root $\al_i$, denote the root vectors 
\begin{equation}
\label{roots}
e_\be = \Ad_{\dot{u}}(e_{\al_i}) \quad \mbox{and} \quad f_\be := \Ad_{\dot{u}}(f_{\al_i}).
\end{equation}

Consider the Schubert cell $B_+ w \cdot B_+$ 
in the full flag variety $G/B_+$ and the isomorphisms
\begin{equation}
\label{ident}
C^-_\ep[w] \cong \KK[U_+ \cap w (U_-) ] \cong \KK [B_+ w \cdot B_+].
\end{equation}
The first one is given by 
\[
f \in \KK [U_+ \cap w (U_-)] \; \; \mt \; \; 
f \big(\exp(z_{\be_1} e_{\be_1}) \ldots \exp(z_{\be_N} e_{\be_N}) \big) \in C^-_\ep[w]
\]
and the second is the pull-back map for the algebraic isomorphism $B_+ w \cdot B_+ \cong U_+ \cap w (U_-)$, 
$g \in U_+ \cap w (U_-) \mt g w \cdot B_+$. (The first isomorphism is the presentation of $U_+ \cap w (U_-)$
as the product of the one-parameter unipotent subgroups of $G$ corresponding to the roots $\be_1, \ldots, \be_N$. The 
coordinate rings of the one-parameter unipotent subgroups are identified with $\KK[z_{\be_j}]$ for $j \in [1,N]$.) 

Denote by $\PP^+$ the set of dominant integral weights of $\g$ and by $\{\vpi_1, \ldots, \vpi_r\}$ the set of fundamental weights. Let
\[
\rho = \vpi_1 + \cdots + \vpi_r.
\]
For $\la \in \PP^+$ and $u, v \in W$, one defines the generalized minors
\[
\De_{u \la, v \la} \in \KK[G] 
\]
as follows. Consider the irreducible highest weight $\g$-module $L(\la)$ with highest weight $\la$. Let $b_\la$ be a highest weight vector 
of $L(\la)$ and $\xi_\la$ be a vector in the dual weight space, normalized by $\lcor \xi_\la, b_\la \rcor =1$. Set
\[
\De_{u \la, v \la} (g) := \lcor \xi_\la, \dot{u}^{-1} g \dot{v} b_\la \rcor, \quad g \in G.
\]
Finally, recall that the support of a Weyl group element $w$ is defined by 
\[
\SS(w) := \{ i \in [1,r] \mid \mbox{$s_i$ occurs in one and thus in any reduced expression of $w$} \}. 
\] 

\bth{Uw-discr} Let $\g$ be a simple Lie algebra, $w$ a Weyl group element and $l>2$ an odd integer
which is $\neq 3$ in the case of $G_2$. Assume that $\KK$ is a field of characteristic 0 which contains a primitive 
$l$-th root of unity $\ep$. Then
\[
d(\UU^-_\ep[w]/C^-_\ep[w]) =_{\KK^\times}  \De_{\rho, w \rho}^L =_{\KK^\times}  \prod_{i \in \SS(w)}\De_{\vpi_i, w \vpi_i}^L
\]  
in the first isomorphism in \eqref{ident} where $L: = l^{N-1}(l-1)$.
\eth
More explicitly, under the first isomorphism in \eqref{ident}, the minor $\De_{\la, w\la}$ corresponds to 
\begin{equation}
\label{prime-e}
\lcor \xi_\la, \exp(z_{\be_1} e_{\be_1}) \ldots \exp(z_{\be_N} e_{\be_N}) \dot{w} b_\la \rcor.
\end{equation}
The equality between the second and third term in \thref{Uw-discr} follows from the product property
\[
\De_{u \la,  v \la} \De_{u \mu, v \la}  = \De_{u(\la+\mu), v(\la+\mu)}, \quad u, v \in W, \la, \mu \in \PP^+
\]
and the fact that $\De_{\vpi_i, w \vpi_i}|_{U_+ \cap w(U_-)} =1$ for $i \notin \SS(w)$.

The algebras $\UU^-_\ep[w]$ and $C^-_\ep[w]$ are defined over $\Qset(\ep)$ and the structure constants for the 
$C^-_\ep[w]$-action on the basis $\YY$ belong to $\Qset(\ep)$ because of \eqref{rel}. This implies that it is sufficient 
to prove \thref{Uw-discr} for any extension $\KK$ of $\Qset(\ep)$. We will do this for $\KK=\Cset$.

{\em{From now on we assume that $\KK= \Cset$ to avoid technicalities with Poisson manifolds over general fields of characteristic 0.}}
(All arguments work for general fields of characteristic 0.)

\thref{Uw-discr} is proved in \S \ref{5.4}. Subsections \ref{5.2} and \ref{5.3} establish results on the 
Poisson geometry of flag varieties and their relation to the induced Poisson structures on $C^-_\ep[w]$. 
Subsection \ref{5.5} applies \thref{Uw-discr} to obtain a formula for the discriminants of the specializations at roots of unity 
of all algebras of rectangular quantum matrices.
%%%%%%%%%%
\subsection{}
\label{5.2} For a $G$-action on a manifold $M$, denote by $\chi \colon \g \to \Ga(M, TM)$ the corresponding infinitesimal 
action and its extension to $\wedge^\bullet \g \to \Ga(M, \wedge^\bullet TM)$.

Let $\De_+$ denote the set of positive roots of $\g$. Recall the definition of the root vectors \eqref{roots} of $\g$.
The standard $r$-matrix for $\g$ is the element
\begin{equation}
\label{stand-r}
r:= \sum_{\be \in \De_+} \frac{\| \be \|^2}{2} e_\be \wedge f_\be \in \wedge^2 \g.
\end{equation}
Define the Poisson bivector field
\[
\pi := - \chi(r) \in \Ga(G/B_+, \wedge^2 T (G/B_+)),
\]
called the standard Poisson structure of the flag variety $G/B_+$. Denote the open Richardson varieties
\[
R_{v,w} := B_- v \cdot B_+ \cap B_+ w \cdot B_+ \subset G/B_+, \quad v \leq w \in W,
\]
see \cite{BL,Deo,Rich}. We will make repeated use of the following facts: 

(A)  The $H$-orbits of symplectic leaves of $(G/B_+, \pi)$ are $R_{v,w}$.

(B) $\ol{R_{v,w}} \cap B_+ w \cdot B_+ = \bigsqcup_{u \in W, u \leq v} R_{u,w}$, 
\\
see \cite[Theorem 4.14]{EL2} and \cite[Theorem 3.2]{Rich}.
\bth{isomP} The composition of the two isomorphisms in \eqref{ident} is an isomorphism of Poisson algebras
\[
(C_\ep^-[w], \{.,.\}) \to (\Cset[B_+ w \cdot B_+], l^2 \ep^{-1} \{.,.\}_\pi).
\]
\eth
For the proof of \thref{isomP} we will need
several constructions for Poisson algebraic groups and Poisson homogeneous spaces, see \cite[Ch. 1]{CP} for background. 
The standard Poisson structure on $G$ is defined by 
\[
\pi_{\st} := \chi_R(r) - \chi_L(r) \in \Ga(G, \wedge^2 TG), 
\]
in terms of \eqref{stand-r}. Here $\chi_R$ and $\chi_L$ denote the infinitesimal actions for the actions of $G$ on itself on the right and the left.
The groups $B_\pm$ are Poisson algebraic subgroups of $(G, \pi_\st)$.
The $r$-matrix for the Drinfeld double of the Poisson algebraic group $(G, \pi_\st)$ is  
\begin{multline*}
r_D:= \sum_{\be \in \De_+} \frac{\| \be \|^2}{2} \left( (e_\be, e_\be) \wedge (f_\be, 0) - (f_\be, f_\be) \wedge (0, e_\be) \right) 
\\
+ \frac{1}{2} \sum_i (h_i, h_i) \wedge (h_i, -h_i) \in (\g \oplus \g)^{\otimes 2}
\end{multline*}
where $\{h_i\}$ is an orthonormal basis of $\Lie(H)$ with respect to the bilinear form $\lcor .,.\rcor$, extending 
the one in \S \ref{5.1}.
The double of $(G, \pi_\st)$ is the group $G \times G$ equipped with the Poisson structure
\[
\pi_D:= \chi_R(r_D) - \chi_L(r_D).
\]
The group 
\[
G^* := \{ (u_- h^{-1}, u_+ h) \mid u_\pm \in U_\pm, h \in H \} \subset G \times G
\]
is a Poisson submanifold of $(G \times G, \pi_D)$; the pair $(G^*, - \pi_D)$ is the 
dual Poisson algebraic group of $(G, \pi_\st)$. The projection onto the first 
component $\eta \colon (G \times G, \pi_D) \to (G, \pi_\st)$,  $\eta (g_1, g_2) = g_1$ is Poisson.
It restricts to the Poisson quotient map
\begin{equation}
\label{eta}
\eta \colon (G^*, -\pi_D) \to (B_-, - \pi_\st).
\end{equation}
 
Denote by $\tau$ and $\theta$ the unique antiautomorphism and automorphism of $G$ which on the Lie 
algebra level are given by 
\[
\tau(e_i) = e_i, \tau(f_i) = f_i, \tau(\al_i\spcheck) = - \al_i\spcheck  \quad \mbox{and} \quad
\theta(e_i) = f_i, \theta(f_i) = e_i, \theta(\al_i\spcheck) = - \al_i\spcheck
\]
for the Chevalley generators of $\g$. It follows from the definition of $\pi_\st$ that 
$\theta \tau \colon (G, \pi_\st) \to (G, - \pi_\st)$ is a Poisson map. This gives rise to the Poisson isomorphism
\begin{equation}
\label{tauth}
\theta \tau \colon (B_-, - \pi_\st) \stackrel{\cong}{\lra} (B_+, \pi_\st). 
\end{equation}
The commutation relations 
\[
\tau \Ad_{\dot{s}_i} = \Ad_{\dot{s}_i^{-1}} \tau \quad \mbox{and} \quad
\theta \Ad_{\dot{s}_i} = \Ad_{\dot{s}_i^{-1}} \theta
\]
and the involutivity of $\tau$ and $\theta$ imply
\begin{equation}
\label{tauth-root}
\theta \tau (f_\be) = e_\be \quad \mbox{for} \; \; \be \in \De_+.
\end{equation}

Recall that an action of a Poisson algebraic group $(G, \pi)$ on a Poisson manifold $(M, \Pi)$ is 
Poisson if the map
\[
(G, \pi) \times (M, \Pi) \to (M, \Pi)
\]
is Poisson. Such a pair $(M, \Pi)$ is a Poisson homogeneous space of $(G, \pi)$ if $M$ is a homogeneous 
$G$-space. If $\Pi$ vanishes at a point $x$, then $(M, \Pi)$ is a Poisson quotient of $(G, \pi)$ via the map
\begin{equation}
\label{qu}
(G, \pi) \to (M, \Pi), \quad g \mt g \cdot x.
\end{equation}
 
The {\em{nonrestricted rational form}} of $\UU_q(\g)$ is the $\Cset[q^{\pm 1}]$-subalgebra,
generated by $E_i$, $F_i$, $K_i^{\pm 1}$ and $(K_i - K_i^{-1})/(q_i - q_i^{-1})$. It will be denoted by $\UU^{\nr}_q(\g)$.
Consider the specialization $\UU_\ep^{\nr}(\g) := \UU_q^{\nr}(\g)/(q - \ep) \UU_q^{\nr}(\g)$ and the canonical projection
$\nu \colon \UU_q^{\nr}(\g) \to \UU_\ep^{\nr}(\g)$. De Concini, Kac and Procesi \cite{DKP1,DKP1b} proved that
$\nu(E_i)^l, \nu(F_i)^l, \nu(K_i)^{\pm l} \in Z(\UU_\ep^{\nr}(\g))$ and that, for good integers $l$, the subalgebra 
$C_\ep^{\nr}(\g)$, generated by them, is a Poisson subalgebra of $Z(\UU_\ep^{\nr}(\g))$ that contains 
all elements $\nu(F_{\be_j})^l$. 

Extend the reduced expression $w = s_{i_1} \ldots s_{i_N}$ to a reduced 
expression $w_\ci= s_{i_1} \ldots s_{i_M}$ of the longest element of $W$ (here 
$M:= \dim \n_-$). Extend the 
set of root vectors $F_{\be_1}, \ldots, F_{\be_N}$ to a set of root vectors $F_{\be_1}, \ldots, F_{\be_N}, \ldots, F_{\be_M}$
by \eqref{root-vect} applied for $j \in [1,M]$. The algebra 
$\UU^-[w_\ci]$ is the $\Cset[q^{\pm 1}]$-subalgebra of $\UU^\nr_q(\g)$ generated by all negative Chevalley generators 
$F_1, \ldots, F_r$. 

By the definition of the induced Poisson structure in \S \ref{3.1}, the embeddings of $\Cset[q^{\pm 1}]$-algebras
$\UU^-[w] \hra \UU^-[w_\ci] \hra \UU_q^{\nr}(\g)$ give rise to the canonical embeddings of Poisson algebras
\begin{equation}
\label{embed2}
(C^-_\ep[w], \{.,.\}) \hra (C^-_\ep[w_\ci], \{.,.\}) \hra
(C^\nr_\ep(\g), \{.,.\})
\end{equation}
where all three Poisson structures are the ones from \eqref{Poisson}.
The first embedding is given by sending $z_{\be_j} \in C^-_\ep[w]$ to $z_{\be_j} \in C^-_\ep[w_\ci]$ for $j \in [1,N]$, recall \eqref{zz}.
The second one is given by $\sig(F_{\be_j})^l \mt \nu(F_{\be_j})^l$.
\medskip

\noindent
{\em{Proof of \thref{isomP}.}} De Concini, Kac and Procesi \cite[Theorem 7.6]{DKP1} constructed an explicit isomorphism of Poisson algebras
\[
I_{DKP} \colon (C^\nr_\ep(\g), \{.,.\}) \stackrel{\cong}{\lra} (\Cset[G^*], -l^2 \ep^{-1} \{.,.\}_{\pi_D}).
\]
It restricts to the Poisson isomorphism 
\[
I_{DKP} \colon (C^-_\ep[w_0], \{.,.\}) \stackrel{\cong}{\lra} (\Cset[F \backslash G^*], -l^2 \ep^{-1} \{.,.\}_{\pi_D})
\]
where $F:= \{ (h^{-1}, h u_+) \mid u_+ \in U_+, h \in H \}$ and $\Cset[F \backslash G^*]$ is viewed
as a Poisson subalgebra of $\Cset[G^*]$. The second isomorphism is explicitly given by $f(z_{\be_j})= z_j$, 
$j \in [1,M]$ where $z_1, \ldots, z_M$ are the coordinate functions on $F \backslash G^*$ from the parametrization
\[
F \backslash G^* = \{ F \cdot \exp(z_M f_{\be_M}) \ldots \exp( z_1 f_{\be_1}) \mid z_1, \ldots, z_M \in \Cset \}.
\]
The explicit statement of this result is given in \cite[Eq. (4.4.1)]{DKP1b}.
The factor $-l^2 \ep^{-1}$ comes from the normalization made in \cite[\S 7.3]{DKP1} and \cite[p. 420]{DKP1b} 
for the induced Poisson bracket on $Z(\UU^\nr_\ep(\g))$. The extra factor of 2 in \cite{DKP1,DKP1b} 
comes from the fact that the Poisson structure $\pi_D$ differs by a factor of 2 from that in \cite{DKP1,DKP1b}. 
Our choice of $\pi_D$ and $\pi_\st$ matches the Poisson structures in \S \ref{4.2}. Composing $I_{DKP}$ with the Poisson maps $\eta^*$ and $\tau^* \theta^*$
(see \eqref{eta} and \eqref{tauth}) gives the Poisson isomorphism
\begin{equation}
\label{composition}
\tau^* \theta^* \eta^* I_{DKP} \colon (C^-_\ep[w_0], \{.,.\}) 
\stackrel{\cong}{\lra}
(\Cset[B_+/H] , l^2 \ep^{-1} \{.,.\}_{\pi_\st})
\end{equation}
where $\Cset[B_+/H]$ is viewed as a Poisson subalgebra of $(\Cset[B_+], l^2 \ep^{-1} \{.,.\}_{\pi_\st})$.
The definition of $I_{DKP}$ and the property \eqref{tauth-root} of $\tau \theta$ imply that the explicit 
form of the isomorphism \eqref{composition} is $\tau^* \theta^* \eta^* I_{DKP}(z_{\be_j}) = \wt{z}_j$, 
$j \in [1,M]$ where $\wt{z}_j$ are the coordinate functions on $\Cset[B_+/H]$ from the parametrization
\[
B_+/H = \{ \exp(\wt{z}_1 e_{\be_1}) \ldots \exp(\wt{z}_M e_{\be_M}) \cdot H \mid z_1, \ldots, z_M \in \Cset \}.
\]

The flag variety $(G/B_+, \pi)$ is a Poisson homogeneous space for $(G, \pi_\st)$. Thus, it is a Poisson 
$(B_+, \pi_\st)$-space. The property (A) in \S \ref{5.2} implies that the Schubert cell 
$(B_+ w \cdot B_+, \pi)$ is a Poisson homogeneous space for Poisson algebraic group $(B_+, \pi_\st)$. By a direct calculation
one checks that $\pi$ vanishes at the base point $w \cdot B_+$. Thus, the fact \eqref{qu} implies that the quotient map 
\[
(B_+, \pi_\st) \to (B_+ w \cdot B_+, \pi), \quad b_+ \mt b_+ w \cdot B_+
\]
is Poisson. In the $\wt{z}_j$ coordinates the map is given by $\exp(\wt{z}_1 e_{\be_1}) \ldots \exp(\wt{z}_M e_{\be_M}) \mt \exp(\wt{z}_1 e_{\be_1}) \ldots \exp(\wt{z}_N e_{\be_N}) w \cdot B_+$. The pull-back map is an embedding of Poisson algebras
\[
(\Cset[B_+ w \cdot B_+], \{.,.\}_\pi) \hra (\Cset[B_+/H], \{.,.\}_{\pi_\st}). 
\]
The theorem follows by combining this embedding, the isomorphism \eqref{composition} and the first embedding in \eqref{embed2}.
\qed
%%%%%%%%%%
\subsection{}
\label{5.3}
Denote by $\QQ$ the root lattice of $\g$.
The algebras $\g$, $\UU_q(\g)$, $\UU^-_\ep[w]$ and $C_\ep^-[w]$ are $\QQ$-graded and the projection $\sig \colon \UU^-[w] \to \UU^-_\ep[w]$ is 
graded. The graded components of these algebras of degree $\ga \in \QQ$ will be denoted by $(.)_\ga$.

\bpr{HomP} The homogeneous prime elements of $(C_\ep^-[w], \{.,.\})$ are
$\De_{\vpi_i, w \vpi_i}$ for $i \in \SS(w)$, in terms of the first identification in \eqref{ident}. They satisfy 
\[
\{\De_{\vpi_i, w \vpi_i}, z \} = - l \ep^{-1} \lcor (w+1) \vpi_i , \ga \rcor \De_{\vpi_i, w \vpi_i} z, \quad
\forall z \in (C_\ep^-[w])_\ga.
\]
\epr
\begin{proof} For $i \in \SS(w)$, the vanishing ideal of $\ol{R_{s_i, w}} \cap B_+ w \cdot B_+$ in $\Cset[B_+ w \cdot B_+]$ is 
$(\De_{\vpi_i, w \vpi_i})$, \cite[Theorem 4.7]{Y-P}. Each of these sets is irreducible and is a union of $H$-orbits of symplectic 
leaves. This follows from the properties (A)-(B) in \S \ref{5.2} and the well known fact that the open Richardson 
varieties $R_{v,w}$ are irreducible. \reref{prime-norm-def} (iii) implies that $\De_{\vpi_i, w \vpi_i} \in C^-_\ep[w]$ 
are homogeneous Poisson prime elements. 

Assume that $f \in C_\ep^-[w]$ is another homogeneous Poisson prime element. 
By \reref{prime-norm-def} (iii), the zero locus $\VV(f)$ of $f$ 
should be a union of $H$-orbits of symplectic leaves of $(B_+ w \cdot B_+, \pi)$. Since 
\[
B_+ w \cdot B_+ = \bigsqcup_{v \in W, v \leq w} R_{v, w} 
\]
and $\dim R_{v,w} = \dim (B_+ w \cdot B_+) - \ell(v)$, either $\VV(f) \cap R_{1,w} \neq \varnothing$ or 
$\VV(f) \cap R_{s_i, w} \neq \varnothing$ for some $i \in \SS(w)$. The first case is impossible since by (A)-(B), $R_{1,w}$ is a 
single $H$-orbit of symplectic leaves and $\ol{R_{1,w}} \supset B_+ w \cdot B_+$. In the second case, 
$\VV(f) \supseteq \ol{R_{s_i, w}} \cap B_+ w \cdot B_+$ because $R_{s_i, w}$ is a single $H$-orbit of leaves. Since $f$ is prime,
$f =_{\Cset^\times} \De_{\vpi_i, w \vpi_i}$.

The formulas for Poisson brackets in the proposition are the specializations at $q=1$ of eq. (5.1) in \cite{Y-IMRN} for $y_1=1$.
\end{proof}
%%%%%%%%%%
\subsection{} 
\label{5.4}
We proceed with the proof of \thref{Uw-discr}. Recall the $C^-_\ep[w]$-basis $\YY$ of $\UU^-_\ep[w]$
from \eqref{Uwbasis}. By \thref{discr-symp} (ii) and \prref{HomP}, 
\begin{equation}
\label{discr-la}
d_{l^N} (\YY : \tr) =_{\Cset^\times} \De_{\la, w \la}
\end{equation}
for some $\la \in \PP^+$. ($C^-_\ep[w]$ is a polynomial algebra and thus a UFD.)
We determine $\la$ by using the methods (1) and (3) in \S \ref{3.3}: We compare the 
degrees of the two sides of the equality (in the $\QQ$-grading) and their Poisson brackets with the elements of $C^-_\ep[w]$.
(Since $\De_{\vpi_i, w \vpi_i}|_{U_+ \cap w(U_-)} =1$ for $i \notin \SS(w)$, $\la$ is only defined up to adding an element of 
$\oplus_{i \notin \SS(w)} \Zset \vpi_i$.) Firstly,
\[
\deg \De_{\la, w \la} = l (w-1) \la.
\]
This follows for instance from \eqref{prime-e} by using that $\deg z_{\be_j} = - l \be_j$. For the reduced expression $w = s_{i_1} \ldots s_{i_N}$, denote
\[
w_{\leq j } := s_{1} \ldots s_{i_j}.
\]
In this notation 
\begin{equation}
\label{rho-be}
- \be_j = - w_{\leq j-1} (\al_{i_j}) = w_{\leq j} \rho - w_{\leq j-1} \rho.
\end{equation}
Since the map $\tr \colon \UU^-_\ep[w] \to C^-_\ep[w]$ is graded,
\begin{align*}
\deg d_{l^N} ( \YY : \tr) &= 2 \sum_{y \in \YY} \deg y = 2 \sum_{k_1, \ldots, k_N =0}^{l-1} 
\deg \sig( F_{\be_1}^{k_1} \ldots F_{\be_N}^{k_N}) 
\\ &= 2 \sum_{k_1, \ldots, k_N =0}^{l-1} k_1 (w_{\leq 1} \rho - \rho)  +  \cdots + k_N (w_{\leq N} \rho - w_{\leq N-1} \rho) =  (l-1) l^N \rho. 
\end{align*}  
Hence, by comparing degrees in \eqref{discr-la}, 
\begin{equation}
\label{ident-1}
(w-1) (\la - (l-1) l^{N-1} \rho) = 0.
\end{equation}
\prref{HomP} and the fact that $\deg z_{\be_j} = - l \be_j$ imply
\begin{equation}
\label{Pbra-De}
\{\De_{\la, w \la}, z_{\be_j} \} = l^2 \ep^{-1} \lcor (w+1) \la , \be_j \rcor \De_{\la, w \la} z_{\be_j}, \quad j \in [1,N].
\end{equation}
To evaluate $\{d_{l^N} (\YY : \tr), z_{\be_j} \}$, we use \prref{Poisson-brack}. Since $\YY$ is a $C^-_\ep[w]$-basis of $\UU^-_\ep[w]$
and $C^-_\ep[w] \cong \Cset[z_{\be_1}, \ldots, z_{\be_N}]$, 
\[
\UU^-_\ep[w] = \oplus_{y \in \YY} \Cset[z_{\be_1}, \ldots, z_{\be_N}] y.
\]
For a monomial $\mu$ in $z_{\be_1}, \ldots, z_{\be_N}$, a basis element $y \in \YY$ and $r \in \UU^-_\ep[w]$, denote by $\coeff_{\mu, y}(r)$ the coefficient 
of $\mu y$ in $r$. For ${\bf{k}}=(k_1, \ldots, k_N) \in \Nset^N$, denote the PBW basis element
\[
F^{\bf{k}} := F_1^{k_1} \ldots F_N^{k_N} \in \UU^-[w]. 
\]

\ble{qSchubert-aux} For all ${\bf{k}} \in [1,l-1]^{\times N}$ and $j \in [1,N]$, 
\[
(\ep^{ \| \al_{i_j} \| /2} - \ep^{-\|\al_{i_j} \| /2})^l \, \coeff_{z_{\be_j}, \sig(F^{\bf{k}})} \big( \del_{F_{\be_j}^l} (\sig(F^{\bf{k}}))  \big) =  
\Big( \sum_{m=1}^N \sign(m-j) k_m \lcor \be_m, \be_j  \rcor \Big) l \ep^{-1}.
\]
\ele
\begin{proof} Consider the right-to-left lexicographic order $\prec$ on $\Nset^N$ given by 
\[
(k_1, \ldots, k_N) \prec (m_1, \ldots, m_N) \quad
\mbox{if $k_N=m_N, \ldots, k_{j+1} = m_{j+1}$ and $k_j < m_j$ for some $j$}.
\]
Recursively applying the straightening law \eqref{rel} gives 
\begin{equation}
\label{filt}
F^{\bf{k}} F^{\bf{m}} = q^{-\sum_{j>a} k_j m_a \lcor \be_j, \be_a \rcor} F^{{\bf{k}}+{\bf{m}}} + \sum_{{\bf{k}}' \prec {\bf{k}} + {\bf{m}}} 
F^{{\bf{k}}'}. 
\end{equation}
Thus,
\[
F_{\be_j}^l F^{\bf{k}} - F^{\bf{k}} F_{\be_j}^l = \big(q^{-l \sum_{a<j} k_j \lcor \be_j, \be_a \rcor} - q^{-l \sum_{a>j} k_j \lcor \be_j, \be_a \rcor}\big) 
F^{{\bf{k}} + l e_j} + \sum_{{\bf{k}}' \prec {\bf{k}} + l e_j} F^{{\bf{k}}'}
\]
where $\{e_1, \ldots, e_N \}$ denotes the standard basis of $\Zset^N$. The lemma follows from this by dividing by $q-\ep$ and applying $\sig$.
\end{proof}
It follows from \eqref{discr-la} and \eqref{Pbra-De} that
\[
\frac{\{d_{l^N} (\YY : \tr), z_{\be_j} \}}{ d_{l^N} (\YY : \tr) z_{\be_j} } \in \Cset.
\]
Now, from \prref{Poisson-brack} we have
\begin{align*}
\frac{\{d_{l^N} (\YY : \tr), z_{\be_j} \}}{ d_{l^N} (\YY : \tr) z_{\be_j} } & = - 2 \sum_{k_1, \ldots, k_N =0}^{l-1} 
(\ep^{ \| \al_{i_j} \| /2} - \ep^{-\|\al_{i_j} \| /2})^l \, \coeff_{z_{\be_j}, \sig(F^{\bf{k}})} \big( \del_{F_{\be_j}^l} (\sig(F^{\bf{k}}))  \big) \\
& = - 2 \sum_{k_1, \ldots, k_N =0}^{l-1} \Big( \sum_{m=1}^N \sign(m-j) k_m \lcor \be_m, \be_j  \rcor \Big) l \ep^{-1} \\
& = (l-1) l^{N+1} \lcor (w+1) \rho, \be_j \rcor \ep^{-1}.
\end{align*}
In the last equality we used the identity $-\lcor(w+1) \rho, \be_j \rcor = \sum_{m=1}^N \sign(m-j) \lcor \be_m, \be_j  \rcor$
which follows from \eqref{rho-be}.
Comparing this with \eqref{Pbra-De}, leads to
\[
\lcor (w+1)(\la - (l-1)l^{N-1} \rho), \be_j \rcor = 0 \quad \mbox{for} \; \; j \in [1,N].
\]
The definition of $\be_j$ implies $\be_j - \al_{i_j} \in \oplus_{m<j} \Zset \al_{i_m}$, and thus,
\[
\bigoplus_{j=1}^N \Zset \be_j = \bigoplus_{i \in \SS(w)} \Zset \al_i.
\] 
Therefore, 
\[
\lcor (w+1)(\la - (l-1)l^{N-1} \rho), \al_i \rcor = 0 \quad \mbox{for} \; \; i \notin \SS(w).
\]
This and the degree formula \eqref{ident-1} give
\[
\lcor \la - (l-1)l^{N-1} \rho, \al_i \rcor = 0 \quad \mbox{for} \; \; i \notin \SS(w),
\]
that is
\[
\la - (l-1) l^{N-1} \rho \in \bigoplus_{i \notin \SS(w)} \Zset \vpi_i.
\]
\thref{Uw-discr} now follows from the fact that $\De_{\vpi_i, w \vpi_i}|_{U_+ \cap w(U_-)} =1$ for $i \notin \SS(w)$.
%%%%%%%%%%%%%
\subsection{}
\label{5.5}
Let $1 \leq m \leq n \in \Zset$. The algebra of (rectangular) quantum matrices $R_q[M_{m, n}]$ is the $\KK[q^{\pm 1}]$-algebra with generators 
$x_{ij}$, $i \in [1,m]$, $j \in [1,n]$ and the four relations in \eqref{qmatr} for $i, k \in [1,n]$, $j,r \in [1,m]$. The algebra is isomorphic to 
$\UU^-[c^m] $ for $\g = {\mathfrak{sl}}_{m+n}$ and the Coxeter element $c = (1 2 \ldots (m+n)) \in S_{m+n}$, see \cite[Proposition 2.1.1]{MC}
and \cite[Lemma 4.1]{Y-IMRN}. The isomorphism $\zeta \colon \UU^-[c^m] \stackrel{\cong}{\lra} R_q[M_{m,n}]$ has the form
\begin{equation}
\label{qisom}
\zeta(F_{\be_k}) = (-q)^{i+j -2} x_{ij}
\end{equation}
for appropriate indices $i \in [1,m]$, $j \in [1,n]$ depending on $k \in [1,mn]$. 
Because of this, the isomorphism in \cite{MC,Y-IMRN}, stated over $\KK(q)$, is valid over $\KK[q^{\pm 1}]$. 

Let $l>2$ be an odd integer and $\ep \in \KK$ be an $l$-th primitive root of unity. Denote the specialization 
$R_\ep[M_{m,n}] := R_q[M_{m, n}]/(q-\ep) R_q[M_{m, n}]$ and the canonical projection $\sig \colon R_q[M_{m, n}] \to R_\ep[M_{m, n}]$.
Analogously to  \S \ref{4.1}, $z_{ij} := \sig(x_{ij})^l \in Z(R_\ep[M_{m,n}])$ and the subalgebra 
$C_\ep[M_{m,n}]$ of $Z(R_\ep[M_{m,n}])$ generated by them is a polynomial algebra in $z_{ij}$. 
The algebra $R_\ep[M_{m,n}]$ is a free module over $C_\ep[M_{m,n}]$ with basis
\[
\{ \sig(x_{11})^{k_{11}} \ldots \sig(x_{mn})^{k_{mn}} \mid 0 \leq k_{11}, \ldots, k_{mn} \leq l-1 \}
\]
where the elements $\sig(x_{ij})$ are listed in any (fixed) order. It follows from \eqref{qisom} 
that $\zeta$ induces in isomorphism $\ol{\zeta} \colon C^-_\ep[c^m] \stackrel{\cong}{\lra} C_\ep[M_{m,n}]$
satisfying
\[
\ol{\zeta}(z_{\be_k}) = (\ep -\ep^{-1})^l (-\ep)^{i+j-2} z_{ij}
\]
where the relation between $i,j$ and $k$ is the same as in \eqref{qisom}. Any minor in these elements is 
a scalar multiple of the corresponding minor in $z_{ij}$. 

Applying \thref{Uw-discr} for $G = {\mathrm{SL}}_{m+n}$ and $w = c^m$, and 
computing explicitly the restriction of the minors $\De_{\vpi_i, c^m \vpi_i}$ to $U_+ \cap c^m(U_-)$, identified with the affine $\KK$-space $M_{m,n}$, 
leads to the following:

\bth{rect-q-matr-discr} Assume that $\KK$ is a field of characteristic 0, $l >2$ is an odd integer, 
$\ep \in \KK$ is a primitive $l$-th root of unity and $m \leq n$ are positive integers. 
The discriminants of the algebras of rectangular quantum matrices are given by
\[
d(R_\ep[M_{m,n}]/C_\ep[M_{m,n}]) =_{\KK^\times} \prod_{j=1}^{m-1} 
\big( \De_{[m-j+1, m]; [1, j]}^L \De_{[1, j] ; [n-j+1, n]}^L \big) 
\prod_{k=1}^{n-m+1} \De_{[1, m]; [k, m+k -1]}^L
\] 
where $L:= l^{mn-1}(l-1)$ and $\De_{I;J}$ denotes the minor in $z_{ij}$ for the rows in the set $I$ and the columns in $J$ as 
in the square case.
\eth
\subsection{}
\label{5.6}
The key properties of a discriminant of an algebra which are used in \cite{BZ,CPWZ1,CPWZ2} are those of being locally dominating or dominating. 
Let $A$ be a finitely generated 
algebra. Fix a generating subspace $V$ of $A$ such that $V \cap \KK = \{0\}$. Denote by $F_j A := ( V \oplus \KK)^j$ the spaces of the induced 
$\Nset$-filtration on $A$. Assume that $\gr A$ is connected and locally finite.
For $a \in A \backslash \{ 0 \}$, denote by $\deg a$ the minimal integer $j$ such that $a \in F_j A$. Fix a $\KK$-basis $a_1, \ldots, a_n$ of $V$
(so, $\deg a_i =1$).

\bde{dom}
An element $f=f(a_1, \ldots, a_n) \in A$ is called locally dominating \cite{CPWZ1} if for every $h \in \Aut(A)$
the following conditions hold:

(i) $\deg f(h(a_1), \ldots, h(a_n)) \geq \deg f$ and 

(ii) if $\deg h(a_i) >1$ for some $i$, then $\deg f(h(a_1), \ldots, h(a_n)) > \deg f$.
\ede

A stronger property of {\em{dominating element}} was defined in \cite{CPWZ1}.

\bde{lead-term} We will say that $f \in A$ has a unique leading term if it has the form
\[
f = t a_1^{k_1} \ldots a_n^{k_n} + \; \; (cwlt), \quad t \in \KK^\times, k_1, \ldots, k_n \in \Zset_+
\]
where the tail $(cwlt)$ is a linear combination of monomials $a_1^{m_1} \ldots a_n^{m_n}$ whose $n$-tuple of powers $(m_1, \ldots, m_n)$ 
is component wise less than $(k_1, \ldots, k_n)$, i.e., $m_i \leq k_i$ for all $i$ and for at least one $i$ the inequality is strict.
\ede
It was proved in \cite[Lemma 2.2]{CPWZ1} that every element that has a unique leading term is dominating and thus locally dominating.
 
The discriminants computed in \cite{BZ,CPWZ1,CPWZ2} posses the stronger property -- those discriminants have a unique 
leading term. The determinants of the matrices $[z_{ij}]_{i,j=1}^n$, considered as elements of $\KK[z_{ij}, 1 \leq i,j \leq n]$,
do not have a unique leading term unless they are of size $1 \times 1$.
It follows from \thref{Uw-discr} that the discriminants of the quantum Schubert cell algebras $d(\UU^-_\ep[w]/C^-_\ep[w]) \in C^-_\ep[w]$ 
do not have a unique leading term unless $\UU^-[w]$ is isomorphic to a quantum affine space algebra.
One can deduce from \cite[Theorem 9.5]{GY-big} that this is precisely the case when $w$ is a subword of a Coxeter element of $W$, (i.e., 
when every simple reflection of $W$ appears at most one time in a reduced expression of $w$). Recently, this subclass of algebras $\UU^-[w]$ 
appeared in a different context (inner faithful Hopf algebra actions) in \cite{EW}.

The next proposition shows that the algebras $\UU^-_\ep[w]$ exhibit a new phenomenon: they can have locally dominating discriminants without possessing 
the stronger property of having a unique leading term. 

\bpr{domina} For all fields $\KK$ of characteristic 0, odd integers integers $l$, and $l$-th roots of unity $\ep \in \KK$, the discriminants
$d(R_\ep[M_2]/C_\ep[M_2])$ and $d(R_\ep[M_{2,3}]/C_\ep[M_{2,3}])$ are locally dominating.
\epr
\begin{proof} Fix an integer $n \geq 2$.
Denote the generators of $R_\ep[M_{2,n}]$ by $a_{ij} := \sig(x_{ij})$. They satisfy the relations \eqref{qmatr} with $q$ replaced by $\ep$. Consider 
the $\Nset$-grading of $R_\ep[M_{2,n}]$ defined by assigning $a_{ij}$ degree one. The grading is connected and locally finite. The grading components 
will be denoted by $R_\ep[M_{2,n}]_k$ and we set $R_\ep[M_{2,n}]_{\geq k} := \oplus_{s \geq k} R_\ep[M_{2,n}]_s$. 
The central elements $z_{ij}$ in \thref{rect-q-matr-discr} are 
$z_{ij} = a_{ij}^l$ and 
\begin{multline}
d(R_\ep[M_{2,n}]/C_\ep[M_{2,n}]) =_{\KK^\times} \\ 
\big( a_{21} ( a_{11} a_{22} - \ep a_{12} a_{21} ) \ldots ( a_{1,n-1} a_{2,n} - \ep a_{1,n} a_{2,n-1})a_{1,n} \big)^{l^{2n}(l-1)}.
\label{De-2f}
\end{multline}
One can show that an identity of this type holds for all $\UU^-[w]$; we give an elementary proof for the case of $R_\ep[M_{2,n}]$. 
Denote $\De:= a_{11} a_{22} - \ep a_{12} a_{21} \in Z(R_\ep[M_2])$. The $q$-binomial formula and
$(a_{11}^{-1} \De)(a_{11}^{-1} a_{12}a_{21}) = \ep^{-2} (a_{11}^{-1} a_{12}a_{21}) (a_{11}^{-1} \De)$ give
\[
a_{22}^l = (a_{11}^{-1} \De + \ep a_{11}^{-1} a_{12}a_{21})^l = (a_{11}^{-1} \De)^l + (a_{11}^{-1} a_{12}a_{21})^l = 
a_{11}^{-l} \big( \De^l + a_{12}^l a_{21}^l \big).
\]
Thus, $\De^l = a_{11}^l a_{22}^l - a_{12}^l a_{21}^l$, and \eqref{De-2f} follows from \thref{rect-q-matr-discr}
and the embeddings \eqref{embed-qm}.

Fix $h \in \Aut(R_\ep[M_{2,n}])$. By \cite[Proposition 3.2]{LL}, it follows from $a_{i1} a_{ij} = \ep a_{ij} a_{i1}$ ($i = 1, 2$, $j\in [2,n]$) that 
$h(a_{ij}) \in R_\ep[M_{2,n}]_{\geq 1}$. Let $h_0(a_{ij})$ be the component of $h(a_{ij})$ in $R_\ep[M_{2,n}]_1$. It follows at once that 
$h_0$ extends to a graded automorphism of $R_\ep[M_{2,n}]$ and 
\[
h(a_{ij}) - h_0(a_{ij}) \in R_\ep[M_{2,n}]_{\geq 2}. 
\]
So,
\begin{align}
\label{propDe}
&h( a_{1,j-1} a_{2,j} - \ep a_{1,j} a_{2,j-1}) -  h_0( a_{1,j-1} a_{2,j} - \ep a_{1,j} a_{2,j-1}) \in R_\ep[M_{2,n}]_{\geq 3},
\\
&h_0( a_{1,j-1} a_{2,j} - \ep a_{1,j} a_{2,j-1}) \in R_\ep[M_{2,n}]_2 \backslash \{ 0 \}.
\nn
\end{align}
for $j \in [1,n-1]$.
As in \deref{dom},  consider the filtration $F_k R_\ep[M_2] := R_\ep[M_2]_{\leq k}$ and, for $a \in R_\ep[M_{2,n}]$, denote 
by $\deg a$ the minimal $k$ such that $a \in R_\ep[M_{2,n}]_{\leq k}$. Eq. \eqref{propDe} implies that 
\begin{equation}
\label{geq-de}
\deg h( a_{1,j-1} a_{2,j} - \ep a_{1,j} a_{2,j-1}) \geq 2, \deg h(a_{ij}) \geq 1 
\end{equation}
for $j \in [1,n-1]$ in the first inequality and $i=1,2$, $j \in [1,n]$ in the second.
This, combined with \eqref{De-2f}, leads to 
\begin{equation}
\label{nonstr-ine}
\deg h( d(R_\ep[M_{2,n}]/C_\ep[M_{2,n}])) \geq 2 n l^{2n} (l-1) = \deg d(R_\ep[M_{2,n}]/C_\ep[M_{2,n}]),
\end{equation}
which verifies condition (i) in \deref{dom} for $d(R_\ep[M_{2,n}]/C_\ep[M_{2,n}])$.
This part of the proof works in the same way for the discriminants of all algebras $\UU^-[w]$ that have the property
that {\em{the $\QQ$-grading of $\UU^-[w]$ can be specialized to an $\Nset$-grading in such way that $\UU^-[w]$ 
is generated in degree 1.}}

Next we verify that $d(R_\ep[M_2]/C_\ep[M_2])$ satisfies the second condition in \deref{dom}. Assume that $h \in \Aut( R_\ep[M_2])$ and $\deg h(a_{ij}) > 1$ 
for at least one generator $a_{ij}$. If $\deg h(a_{21}) >1$ or $\deg h(a_{12})>1$, then \eqref{geq-de}
implies that the inequality \eqref{nonstr-ine} is strict. Assume that $\deg h(a_{12}) =\deg h(a_{21})=1$. Then either $\deg h(a_{11}) >1$ 
or $\deg h(a_{22}) >1$, so, $\deg h(a_{11} a_{22}) >2$ and $\deg h (a_{12} a_{21}) =2$. Thus,
\[
\deg h(a_{11} a_{22} - \ep a_{12} a_{21}) >2,
\]
which implies that the inequality \eqref{nonstr-ine} is strict. This proves that $d(R_\ep[M_2]/C_\ep[M_2])$ satisfies 
the condition (ii) in \deref{dom}.

Finally, let $h \in \Aut ( R_\ep[M_{2,3}])$ and $\deg h(a_{ij}) > 1$ for at least one generator $a_{ij}$. 
If $\deg h(a_{13}) >1$ or $\deg h(a_{21})>1$, then, again \eqref{geq-de} implies that the inequality \eqref{nonstr-ine} is strict.
Consider the case $\deg h(a_{13}) = 1$ and $\deg h(a_{21})=1$; so, $h(a_{13}) = h_0(a_{13})$ and $h(a_{21})= h_0(a_{21})$.
We prove by contradiction that $d(R_\ep[M_{2,3}]/C_\ep[M_{2,3}])$ satisfies 
the condition (ii) in \deref{dom}. Assume that this is not the case, then
$\deg h( d(R_\ep[M_{2,3}]/C_\ep[M_{2,3}])) = 6 l^6 (l-1)$ and
\begin{align}
\label{h-2id}
&h(a_{11}) h (a_{22}) - \ep h (a_{12}) h(a_{21}) \in R_\ep[M_{2,3}]_2, 
\\
&h(a_{12}) h (a_{23}) - \ep h (a_{13}) h(a_{22}) =
h (a_{23}) h(a_{12}) - \ep^{-1} h (a_{13}) h(a_{22}) \in R_\ep[M_{2,3}]_2.
\nn
\end{align}
Therefore, 
\begin{align*}
&\deg h(a_{11}) + \deg h (a_{22}) = 1 + \deg h (a_{12}), \\
&\deg h(a_{12}) + \deg h (a_{23}) =  1+  \deg h(a_{22}),
\end{align*} 
which imply that $\deg h(a_{11}) = \deg h(a_{23}) =1$ and $\deg h (a_{12}) = \deg h(a_{22}) =k > 1$ (the latter holds because 
of the assumption that $\deg h(a_{ij}) > 1$ for at least one generator $a_{ij}$). Denote by $b_{12}$ and $b_{22}$ the (nonzero) components 
of $h(a_{12})$ and $h(a_{22})$ in $R_\ep[M_{2,3}]_k$. From \eqref{h-2id} we obtain
\[
h_0(a_{11}) b_{22} - \ep h_0 (a_{21}) b_{12} = 0 = h_0 (a_{23}) b_{12}  - \ep^{-1} h_0 (a_{13}) b_{22}.
\]
Thus, $b_{22} b_{12}^{-1} = \ep h_0(a_{11})^{-1} h_0(a_{21}) = \ep h_0(a_{13})^{-1} h_0(a_{23})$ and 
$a_{11}^{-1} a_{21} = a_{13}^{-1} a_{23}$ which is a contradiction. Therefore, $d(R_\ep[M_{2,3}]/C_\ep[M_{2,3}])$ satisfies 
the condition (ii) in \deref{dom}.
\end{proof}
%%%%%%%%%%%%%%%%%%%%%% References %%%%%%%%%%%%%%%%%%%%%%%%%%%%%%%%%%%%%%%

%%%%%%%%%%%%%%%%%%%%%%%%%%%%%%%%%%%%%%%%%%%%%%%%%%%%%%%%%%%%%%%%%%%%%%%%%%%%%%%
%%%%%%%%%%%%%%%%%%%%%%%%%%%%%%%%%%%%%%%%%%%%%%%%%%%%%%%%%%%%%%%%%%%%%%%%%%%%%%
\end{document}